\documentclass{amsart}
\usepackage{amsmath,amssymb,amsthm}
\usepackage{mathabx,bbm,esint}
\usepackage{xcolor}
\usepackage{url,hyperref}
\hypersetup{colorlinks=false,pdfborder={0 0 0.6},
 linkbordercolor={1 0 0},citebordercolor={0 1 0},
 urlbordercolor={0 1 0},filebordercolor={0 1 0},
 pdftitle={Restriction Estimates Using Multilinear Kakeya},
 pdfauthor={Xiangyu Wang},bookmarksnumbered=true}

\numberwithin{equation}{section}
\newtheorem{theorem}{Theorem}[section]
\newtheorem{lemma}[theorem]{Lemma}
\newtheorem{proposition}[theorem]{Proposition}
\theoremstyle{definition}
\newtheorem{definition}[theorem]{Definition}
\theoremstyle{remark}

\newcommand{\R}{\mathbb R}
\newcommand{\T}{\mathbb T}
\newcommand{\1}{\mathbbm 1}
\newcommand{\BL}{\mathrm{BL}}
\newcommand{\RD}{\operatorname{RapDec}}
\newcommand{\supp}{\operatorname{supp}}
\newcommand{\ang}{\measuredangle}
\newcommand{\eps}{\epsilon}
\newcommand{\cB}{\mathcal B}
\newcommand{\cQ}{\mathcal Q}
\newcommand{\cG}{\mathcal G}
\newcommand{\step}[1]{\medskip\noindent\textbf{#1}}
\newcommand{\orcid}[1]{\href{https://orcid.org/#1}{\texttt{ORCID: #1}}}
\begin{document}

\title[Restriction estimates using multilinear Kakeya]{Restriction Estimates Using Multilinear Kakeya}

\author{Xiangyu Wang} \address{Department of Mathematics, University of Illinois Urbana-Champaign\newline \orcid{0009-0003-5983-6961}} \email{xw70@illinois.edu}

\begin{abstract}
We use multilinear Kakeya estimates to establish restriction estimates for the paraboloid in dimensions $n\ge5$. The proof inserts Zahl's multilinear Kakeya estimate for direction-separated tubes into the two-ends method of Wang and Wu by means of a broad--narrow argument. In several dimensions, including $n=9$, the resulting exponents were the best known before the recent proof of the local smoothing conjecture by OpenAI.
\end{abstract}
\maketitle


\section{Introduction}\label{sec:introduction}
Let $n\ge2$, let $B^{n-1}$ be the unit ball in $\R^{n-1}$, and define the extension operator for the paraboloid by
\[
 Ef(x',x_n)=\int_{B^{n-1}}e^{i(x'\cdot\xi+x_n|\xi|^2)}f(\xi)\,d\xi,
 \qquad (x',x_n)\in\R^{n-1}\times\R.
\]
The restriction problem asks for the range of $p$ for which
\begin{equation}\label{eq:diagonal}
 \|Ef\|_{L^p(\R^n)}\le C_{n,p}\|f\|_{L^p(B^{n-1})}.
\end{equation}
The conjectured range is $p>2n/(n-1)$. Polynomial partitioning has led to a sequence of advances in this problem; see, for example, \cite{Guth,HZ,GWZ}. Wang and Wu \cite{WW} introduced a different approach, based on refined decoupling and on incidence estimates for tubes under a \emph{two-ends} condition. In this paper we combine their approach with the broad--narrow method of Bourgain and Guth \cite{BG,Guth} and with Zahl's multilinear Kakeya estimate for direction-separated tubes \cite{Zahl}.

To state the result, let $n\ge5$ and fix
\begin{equation}\label{eq:k-range}
 2\le k\le\left\lfloor\frac{n+2}{2}\right\rfloor,
 \qquad
 d=d_{n,k}:=\frac{n^2+k^2+n-k}{2n},
\end{equation}
and put $D=(4n-5)(d-1)+(n-1)(n-3)$. Define
\[
 p_b(n,k)=2+\frac{4(2d+n-5)}D,\qquad p_N(n,k)=2+\frac4{2n-k}.
\]
The number $d_{n,k}$ is the exponent in Zahl's $k$-linear Kakeya estimate (Lemma~\ref{thm:zahl-input}). The exponent $p_b$ is the threshold for the $k$-broad estimate that we prove in Section~\ref{sec:induction}, and $p_N$ is the threshold above which a $k$-broad estimate implies a linear one (Section~\ref{sec:narrow}).

\begin{theorem}\label{thm:main}
For $n\ge5$, the estimate
\eqref{eq:diagonal} holds whenever
\[
 p>\min_{2\le k\le\lfloor(n+2)/2\rfloor}
 \max\{p_b(n,k),p_N(n,k)\}.
\]
\end{theorem}

Table~\ref{tab:thresholds} lists the exponent of Theorem~\ref{thm:main} for $5\le n\le12$. For large $n$ the minimum is attained when $k/n$ is close to the real root $x_0=0.4534\ldots$ of $x^3+2x-1=0$, and the exponent is $2+\frac{4}{(2-x_0)n}+O(n^{-2})=2+\frac{2.58631445\ldots}{n}+O(n^{-2})$.

\begin{table}[ht]
\caption{The exponent of Theorem~\ref{thm:main} and a value of $k$ that attains the minimum.}\label{tab:thresholds}
\small\setlength{\tabcolsep}{3.5pt}
\begin{tabular}{c|cccccccc}
$n$ & $5$ & $6$ & $7$ & $8$ & $9$ & $10$ & $11$ & $12$\\\hline
$k$ & $3$ & $3$ & $4$ & $4$ & $5$ & $5$ & $6$ & $6$\rule{0pt}{2.6ex}\\[0.3ex]
exponent & $2+\frac{144}{235}$ & $2+\frac12$ & $2+\frac{328}{789}$ & $2+\frac{216}{599}$ & $2+\frac{292}{929}$ & $2+\frac{144}{511}$ & $2+\frac{24}{95}$ & $2+\frac{40}{173}$\\[0.6ex]
 & $2.6128$ & $2.5000$ & $2.4157$ & $2.3606$ & $2.3143$ & $2.2818$ & $2.2526$ & $2.2312$
\end{tabular}
\end{table}

\subsection{Relation to other results}\label{subsec:relation}
Recently, OpenAI released the manuscript \cite{OAI}, which gives a proof of Sogge's local smoothing conjecture for the Euclidean wave equation in every dimension. Local smoothing implies the Bochner--Riesz conjecture, which implies the restriction conjecture \cite{Tao99}; accordingly, \cite{OAI} obtains \eqref{eq:diagonal} for all $p>2n/(n-1)$, which supersedes Theorem~\ref{thm:main}. The present paper was made public before \cite{OAI}, and we believe that its method, which combines the two-ends approach of Wang and Wu with multilinear Kakeya estimates, has independent mathematical interest. Below we compare Theorem~\ref{thm:main} with the results that were available before \cite{OAI}.

For $n=9$, $11$ and $13$, the exponents of Theorem~\ref{thm:main} are $2.3143\ldots$, $2.2526\ldots$ and $2.2111\ldots$, whereas those recorded in \cite[Figure~2]{HZ} are $2.3166\ldots$, $2.2536\ldots$ and $2.2115\ldots$. To the best of our knowledge, Theorem~\ref{thm:main} gave the best exponent among the results published before \cite{OAI} in these dimensions, and also in several others, for example $n=5,7,8,10$. We point out, however, that the broad estimate \cite[Theorem~4.2]{GWZ} of Guo, Wang and Zhang is stated only for $2n/5\le k\le n/2$; a numerical optimization of the parameters in its proof, for all $k$, gives better exponents than Theorem~\ref{thm:main} for $n=11$ and $n=13$, though not for $n=9$. Finally, Wang and Wu \cite[Theorem~0.17]{WW} obtained \eqref{eq:diagonal} for $p>\frac{154n+6}{77n-95}=2+\frac{28}{11n}+O(n^{-2})$, using a two-ends inequality of Katz and Tao \cite{KT}. This exponent is better than that of Theorem~\ref{thm:main} for $n=28$ and $n\ge30$, and worse for $5\le n\le27$ and $n=29$.

\subsection{Notation and conventions}\label{subsec:notation}
We write $A\lesssim B$ when $A\le CB$ with a constant $C$ that may depend on $n$, $k$, $p$ and on the parameters $\eta,\kappa,a_0,\sigma,\alpha,b_0,\eps_0,\eps$ fixed in the proofs, but not on $R$, $r$, $\delta$, $K$, $A$ or on the functions and families of tubes under consideration; $A\sim B$ means $A\lesssim B\lesssim A$. The notation $A\lessapprox B$ is used only in Section~\ref{subsec:strategy}, where it means that $A\le C_\eps r^{\eps}B$ for every $\eps>0$; in the rest of the paper, all losses of the form $r^{\eps}$, $r^{O(\alpha)}$, $K^{O(1)}$ and $(\log r)^{O(1)}$ are written explicitly.

A term $\RD(r)$ is a quantity bounded by $C_Lr^{-L}$ for every $L\ge1$. We write $|S|$ for the Lebesgue measure and $\#S$ for the cardinality of a set $S$, $N_\rho(S)$ for its $\rho$-neighborhood, and $B_r$ for a ball of radius $r$ in $\R^n$. A \emph{cap} is a ball in $\R^{n-1}$ of radius at most $1$. For vectors $v_1,\ldots,v_k\in\R^n$, $|v_1\wedge\cdots\wedge v_k|$ is the $k$-dimensional volume of the parallelepiped that they span.

We use dyadic pigeonholing in the following form. If a finite set $S$ is partitioned into $L$ classes and $w\colon S\to[0,\infty)$, then some class $S'$ satisfies $\sum_{S'}w\ge L^{-1}\sum_Sw$. In our applications the classes are defined by the dyadic values of one or several quantities, each of which takes values in an interval $[r^{-C},r^{C}]$, so that $L\lesssim(\log r)^{C}$.

\subsection{Strategy of the proof}\label{subsec:strategy}
We first recall the method of Wang and Wu \cite{WW} in a simplified form. Decompose $f$ into wave packets $f_T$ at scale $r\ge1$: each $f_T$ is supported in a cap of radius $r^{-1/2}$, and inside a ball $B_r$ of radius $r$ the function $Ef_T$ is concentrated on a tube $T$ of length $r$ and radius $r^{1/2}$. Let $X\subset B_r$ be a union of cubes of side $r^{1/2}$, and suppose that the norm $\|f_T\|_2$, the number $m$ of wave packets in each direction, the number $\mu$ of tubes that meet a cube of $X$, and the \emph{density} $\lambda=|T\cap X|/|T|$ are essentially constant. We have the elementary estimate and the refined decoupling estimate of Guth, Iosevich, Ou and Wang \cite{GIOW}:
\begin{equation}\label{eq:intro-L2}
 \int_X|Ef|^2\lesssim\lambda r\,\|f\|_2^2,\qquad
 \int_X|Ef|^q\lessapprox r^{-1}\Bigl(\frac\mu m\Bigr)^{\frac2{n-1}}\|f\|_2^2\,M_r(f)^{q-2},
\end{equation}
where $q=\frac{2(n+1)}{n-1}$ and $M_r(f)$ is the largest $L^2$-average of $f$ over a cap of radius $r^{-1/2}$. To combine them one needs an \emph{incidence estimate}
\begin{equation}\label{eq:intro-ab}
 \frac\mu m\lessapprox\delta^{-a}\lambda^{-b},\qquad \delta=r^{-1/2}.
\end{equation}
If \eqref{eq:intro-ab} is inserted into the second estimate of \eqref{eq:intro-L2}, and the result is interpolated with the first estimate so that the powers of $r$ cancel, one obtains $p=2+\frac4{2(n-1)-a}$, and the bound for $\int_X|Ef|^p$ carries the factor $\lambda^{c(n-1-a-2b)}$ with $c>0$. Thus \eqref{eq:intro-ab} suffices by itself if $a+2b\le n-1$, and the exponent then depends only on $a$; if $a+2b>n-1$, the bound deteriorates as $\lambda\to0$, and \eqref{eq:intro-ab} has to be combined with another estimate. For general sets the loss in $\lambda$ is large: a Kakeya maximal estimate at dimension $d'$ gives $(a,b)=(n-d',d'-1)$. The observation of Wang and Wu is that an induction on scales allows one to assume that $X\cap T$ is not concentrated on a short segment of $T$ (the two-ends condition, Definition~\ref{def:twoends}), and that much smaller values of $b$ are then available.

We use two incidence estimates.
\begin{itemize}
\item[(i)] \emph{The two-ends hairbrush estimate} (Lemma~\ref{thm:GWX-input}). Wolff's hairbrush argument \cite{Wolff} groups the tubes that meet a fixed tube according to the $2$-planes that they span with it, estimates each group in the plane, and gives $(a,b)=\bigl(\frac{n-2}2,\frac n2\bigr)$. Wang and Wu replaced the planar estimate by their two-ends Furstenberg inequality \cite[Theorem~2.1]{WW}, which loses only $\lambda^{-1/2}$; this gives $(a,b)=\bigl(\frac12,\frac34\bigr)$ in $\R^3$ \cite[Proposition~3.2]{WW}, and Gao, Wu and Xi \cite[Lemma~4.1]{GWX} obtained $(a,b)=\bigl(\frac{n-2}2,\frac34\bigr)$ in every dimension. Here $a+2b=\frac{n+1}2\le n-1$, so this estimate suffices by itself and gives $p>2+\frac8{3n-2}$; for $n>3$ it leaves a positive power of $\lambda$ unused.
\item[(ii)] \emph{Zahl's multilinear Kakeya estimate} (Lemma~\ref{thm:zahl-input}). Zahl \cite{Zahl} proved that direction-separated tubes cannot cluster near algebraic varieties of low degree, and deduced that the $k$-linear Kakeya theorem of Bennett, Carbery and Tao \cite{BCT} and Guth \cite{GuthEndpoint} holds at dimension $d_{n,k}=k+\frac{(n-k)(n-k+1)}{2n}$, instead of $k$, for such tubes. For $k$ transverse families it gives \eqref{eq:intro-ab} with $(a,b)=(n-d,d-1)$, where $\mu$ and $\lambda$ are geometric means over the families (Lemma~\ref{lem:geometry}). Its value of $a$ is smaller than in (i) when $d>\frac{n+2}2$, that is, when $k(k-1)>n$. But it is an estimate for complete tubes and does not use the two-ends condition, and $a+2b=n+d-2>n-1$, so it is too weak for sparse sets.
\end{itemize}
This is why the deeper estimate (ii) does not make (i) redundant. The argument has to handle all densities $\lambda$ at once, and only (i) controls sparse sets; (ii) is used to improve the power of $r$ at the expense of the unused power of $\lambda$ in (i). We apply the refined decoupling estimate once with (i) and once with (ii), and combine the two resulting bounds with the $L^2$ estimate, with weights chosen so that the powers of $r$ and of $\lambda$ both cancel (see \eqref{multilinear-norm} and the lines that follow it). This gives the exponent $p_b(n,k)$.

Since (ii) requires transversality, we work with the $k$-broad norms of Guth \cite{Guth} (Section~\ref{sec:broad}), which are controlled by the contributions of $k$ transverse caps (Lemma~\ref{eq:extract-bound}). The main step is the $k$-broad estimate of Proposition~\ref{prop:broad}, which we prove by the two-ends induction on scales of \cite{WW}. The passage to \eqref{eq:diagonal} is the broad--narrow argument of Bourgain and Guth \cite{BG,Guth}: it requires $p>p_N(n,k)$, and the decoupling inequality that it uses requires $p\le\frac{2(k-1)}{k-2}$, which is why $k$ is restricted as in \eqref{eq:k-range}. We use Zahl's multilinear estimate directly in the broad argument, rather than the linear Kakeya maximal estimate that he deduces from it \cite[Lemma~1.4]{Zahl}, which requires $d\le n-k+2$; this avoids the loss that can occur in that conversion.

Section~\ref{sec:preliminaries} collects definitions and tools, Section~\ref{sec:newframework} deduces from Zahl's theorem the incidence estimate for transverse families, and Section~\ref{sec:proof} contains the proof of Theorem~\ref{thm:main}.

\subsection*{Acknowledgements} Claude was employed to refine the English writing and to help check the details of the proof. The author thanks an anonymous referee for pointing out that, in an earlier version of this manuscript, the author's uncritical reliance on ChatGPT for the writing had made the text very difficult to read.

\section{Preliminaries}\label{sec:preliminaries}
\subsection{Tubes, shadings, and the two-ends condition}\label{subsec:tubes}
Throughout this subsection and the next, $0<\delta<1$.

\begin{definition}[Tubes and directions]\label{def:tubes}
A \emph{$\delta$-tube} is a set $T\subset\R^n$ such that
\[
 N_{\delta/C_n}(\ell)\subset T\subset N_{C_n\delta}(\ell)
\]
for some line segment $\ell=\ell_T$ whose length lies between $1$ and $C_n$. Here $C_n\ge1$ is a constant that depends only on $n$ and is fixed throughout the paper. The \emph{direction} $v(T)\in S^{n-1}$ of $T$ is the unit vector parallel to $\ell_T$ whose last coordinate is positive; all tubes in this paper satisfy $v(T)_n\ge1/5$.
\end{definition}

The constant $C_n$ is needed because the tubes produced by the wave packet decomposition in Section~\ref{sec:packets} are sheared cylinders, and not exact neighborhoods of segments. The estimates quoted in Section~\ref{subsec:geom-tools} are stated in the literature for $\delta$-neighborhoods of unit segments; they extend to $\delta$-tubes in the sense of Definition~\ref{def:tubes} after the tubes are enlarged and rescaled by factors that depend only on $n$, and this affects only the constants.

To count tubes with nearly the same direction, we fix a maximal $\delta$-separated subset $\Omega_\delta$ of $\{\omega\in S^{n-1}:\omega_n\ge1/5\}$, and for each $\delta$-tube $T$ we fix a point $\omega(T)\in\Omega_\delta$ with $|v(T)-\omega(T)|<\delta$.

\begin{definition}[Parallel multiplicity and direction separation]\label{def:parallel}
Let $\T$ be a finite family of $\delta$-tubes.
\begin{itemize}
\item[(i)] For $\omega\in\Omega_\delta$, the set $\{T\in\T:\omega(T)=\omega\}$ is the \emph{direction class} of $\omega$ in $\T$. For an integer $m\ge1$, the family $\T$ is \emph{$m$-parallel} if each of its direction classes contains at most $m$ tubes.
\item[(ii)] The family $\T$ is \emph{direction-separated} if $|v(T)-v(T')|\ge\delta$ for all distinct $T,T'\in\T$.
\end{itemize}
\end{definition}

An $m$-parallel family satisfies
\begin{equation}\label{eq:card-parallel}
 \#\T\le m\,\#\Omega_\delta\lesssim m\delta^{-(n-1)}.
\end{equation}
A $1$-parallel family need not be direction-separated, because two tubes in neighboring direction classes may be almost parallel. This is repaired by splitting $\Omega_\delta$ into a bounded number of subsets. Since $\Omega_\delta$ is $\delta$-separated, each of its points has at most $c_n-1$ other points of $\Omega_\delta$ within distance $3\delta$, where $c_n$ depends only on $n$; hence there is a partition
\begin{equation}\label{eq:coloring}
 \Omega_\delta=\Omega_\delta^{(1)}\cup\cdots\cup\Omega_\delta^{(c_n)}
 \quad\text{with}\quad
 |\omega-\omega'|\ge3\delta\ \text{ for all distinct }\omega,\omega'\in\Omega_\delta^{(j)},\ 1\le j\le c_n.
\end{equation}
If $\T$ is $1$-parallel and $T,T'\in\T$ are distinct tubes with $\omega(T),\omega(T')\in\Omega_\delta^{(j)}$, then $|v(T)-v(T')|>|\omega(T)-\omega(T')|-2\delta\ge\delta$. Therefore
\begin{equation}\label{eq:one-parallel}
 \text{every $1$-parallel family is the union of $c_n$ direction-separated families.}
\end{equation}

\begin{definition}[Shadings]\label{def:shading}
Let $\cG_\delta$ be the grid of cubes $\delta(z+[0,1)^n)$, $z\in\mathbb Z^n$. A \emph{shading} of a family $\T$ of $\delta$-tubes is a map $Y$ that assigns to each $T\in\T$ a set $Y(T)$ which is a union of cubes of $\cG_\delta$ that meet $T$. We write
\[
 E_Y=\bigcup_{T\in\T}Y(T),\qquad \mu_Y(x)=\sum_{T\in\T}\1_{Y(T)}(x).
\]
The number $\mu_Y(x)$ is the \emph{multiplicity} of the shading at $x$; it is constant on each cube of $\cG_\delta$. For $0<\lambda\le1$, the shading $Y$ is \emph{$\lambda$-dense} if $|Y(T)|\ge\lambda|T|$ for every $T\in\T$. The \emph{average density} of $Y$ is the number $\bar\lambda$ defined by $\sum_{T\in\T}|Y(T)|=\bar\lambda\sum_{T\in\T}|T|$.
\end{definition}

\begin{definition}[Two ends]\label{def:twoends}
Let $0<\eps_2<\eps_1<1$. A shading $Y$ of $\T$ is \emph{$(\eps_1,\eps_2)$-two-ends} if
\[
 |Y(T)\cap B(x,\delta^{\eps_1})|
 \le \delta^{\eps_2}|Y(T)|
\]
for every $T\in\T$ and every $x\in\R^n$.
\end{definition}

Thus a shading is two-ends if no ball of radius $\delta^{\eps_1}$, which is much shorter than the tube, contains more than the fraction $\delta^{\eps_2}$ of $Y(T)$. The definitions in \cite[Definition~1.21]{WW} and \cite[Definition~3.6]{GWX} use segments of $T$ of length $\delta^{\eps_1}$ in place of balls; such a segment is contained in a ball of radius $\delta^{\eps_1}$, so a shading that is two-ends in the sense of Definition~\ref{def:twoends} is two-ends in their sense.

\subsection{Two geometric estimates}\label{subsec:geom-tools}
The first estimate bounds the multiplicity of a two-ends shading.

\begin{lemma}[Two-ends hairbrush estimate; Gao--Wu--Xi]\label{thm:GWX-input}
Let $n\ge3$, $0<\eps_2<\eps_1<1$ and $\eps>0$. Let $\T$ be an $m$-parallel family of $\delta$-tubes in $\R^n$, and let $Y$ be a $\lambda$-dense, $(\eps_1,\eps_2)$-two-ends shading of $\T$. Then there is a set $E_H\subset E_Y$, which is a union of cubes of $\cG_\delta$, such that
\[
|E_Y\setminus E_H|\lesssim\delta^{\epsilon_1}|E_Y|,
\qquad
\mu_Y(x)\le C\delta^{-\eps-C_0\eps_1}\,
m\lambda^{-3/4}\delta^{-(n-2)/2}
\quad\text{for }x\in E_H.
\]
Here $C_0$ depends only on $n$, and $C$ and the implicit constant depend only on $n,\eps,\eps_1,\eps_2$.
\end{lemma}

This is a restatement of \cite[Lemma~4.1]{GWX}, and we indicate the changes. The lemma in \cite{GWX} is formulated for families that obey Wolff's axiom with error $m$ \cite[Definition~2.3]{GWX}: every rectangular box of dimensions $1\times s_1\times\cdots\times s_{n-1}$ contains $\lesssim m\,s_1\cdots s_{n-1}\delta^{-(n-1)}$ tubes of the family. An $m$-parallel family has this property, because the directions of the tubes contained in such a box belong to $\lesssim s_1\cdots s_{n-1}\delta^{-(n-1)}$ direction classes. The lemma in \cite{GWX} gives the multiplicity bound
\[
 \mu_Y(x)\le C\delta^{-\eps-C_0\eps_1}\,m^{\frac n{2n-2}}\lambda^{-3/4}\delta^{-(n-2)/2}\,(\delta^{n-1}\#\T)^{\frac{n-2}{2n-2}}
\]
outside a subset of $E_Y$ of measure at most $\delta^{\eps_1}|E_Y|$, and the bound in Lemma~\ref{thm:GWX-input} follows from \eqref{eq:card-parallel}. If $\mu_0$ denotes the right-hand side of the bound in Lemma~\ref{thm:GWX-input}, we may take $E_H=\{x\in E_Y:\mu_Y(x)\le\mu_0\}$; this set is a union of cubes of $\cG_\delta$, because $\mu_Y$ is constant on these cubes. As explained in Section~\ref{subsec:strategy}, the lemma is a refinement of Wolff's hairbrush bound \cite{Wolff} for two-ends shadings; the case $n=3$ is \cite[Proposition~3.2]{WW}.

The second estimate is Zahl's multilinear Kakeya theorem \cite[Theorem~1.3]{Zahl}. Recall that $d_{n,k}$ was defined in \eqref{eq:k-range}.

\begin{lemma}[Multilinear Kakeya for direction-separated tubes; Zahl]
\label{thm:zahl-input}
Let $2\le k\le n$, $d=d_{n,k}$ and $\eps>0$. For every direction-separated family $\T$ of $\delta$-tubes in $\R^n$,
\[
\begin{aligned}
 &\left\|\sum_{T_1,\ldots,T_k\in\T}
 \1_{T_1}\cdots\1_{T_k}
 |v(T_1)\wedge\cdots\wedge v(T_k)|^{k/d}
 \right\|_{L^{d/(k(d-1))}(\R^n)}^{1/k}\\
 &\qquad\qquad\le C_{n,\eps}\,
 \delta^{1-n/d-\eps}
 (\delta^{n-1}\#\T)^{\frac{n(d-1)}{(n-1)d}}.
\end{aligned}
\]
\end{lemma}

The left-hand side counts the $k$-tuples of tubes through a point, with a weight that vanishes when the $k$ directions are linearly dependent. For comparison, the Kakeya maximal estimate at dimension $d$ is the bound $\|\sum_{T\in\T}\1_T\|_{L^{d/(d-1)}}\le C_\eps\delta^{1-n/d-\eps}(\delta^{n-1}\#\T)^{(d-1)/d}$ for direction-separated families. Thus Lemma~\ref{thm:zahl-input} is an estimate of this type at dimension $d_{n,k}$ for the part of the sum that comes from $k$ tubes with independent directions. The form in which we use it is Lemma~\ref{lem:geometry} below.

\subsection{Broad norms}\label{sec:broad}
We recall the $k$-broad norms of Guth \cite{Guth}. Fix $R\gg1$ and a dyadic number $K\sim R^\kappa$, where $\kappa>0$ is small.

\emph{Caps and families.} Let $\tau$ range over a finitely overlapping cover of $B^{n-1}$ by caps of radius $c_{n,k}K^{-1}$ with centers $\xi_\tau\in B^{n-1}$, where $c_{n,k}>0$ is a sufficiently small constant, fixed in Lemma~\ref{eq:extract-bound}. We still call them $K^{-1}$-caps, and $2\tau$ and $3\tau$ denote the caps with the same center and two or three times the radius. Choose a smooth partition of unity $\sum_\tau\psi_\tau=1$ on $B^{n-1}$ with $0\le\psi_\tau\le1$ and $\supp\psi_\tau\subset\tau$, and for $f\in L^2(B^{n-1})$ put $f_\tau=f\psi_\tau$.

The wave packet decompositions below are applied to each $f_\tau$ separately, and they enlarge the supports slightly. For this reason we work with \emph{families} $f=(f_\tau)_\tau$ of functions $f_\tau\in L^2(\R^{n-1})$ with $\supp f_\tau\subset2\tau$, indexed by the caps. For such a family,
\[
 Ef_\tau(x',x_n)=\int_{\R^{n-1}}e^{i(x'\cdot\xi+x_n|\xi|^2)}f_\tau(\xi)\,d\xi,
 \qquad
 \|f\|_2^2:=\sum_\tau\|f_\tau\|_2^2 .
\]
Sums and differences of families are formed cap by cap. A function $f\in L^2(B^{n-1})$ is identified with the family $(f\psi_\tau)_\tau$; for this family, $\sum_\tau\|f\psi_\tau\|_2^2\sim\|f\|_{L^2(B^{n-1})}^2$.

\emph{Directions.} For $\xi\in\R^{n-1}$ let
\[
 G(\xi)=\frac{(-2\xi,1)}{|(-2\xi,1)|}\in S^{n-1}.
\]
This is the direction of the tubes on which the wave packets with frequencies near $\xi$ are concentrated; see Section~\ref{sec:packets}. For a set $S\subset\R^{n-1}$ we write $G(S)=\{G(\xi):\xi\in S\}$, and for a linear subspace $V\subset\R^n$ we write $\ang(G(S),V)$ for the smallest angle between a vector in $G(S)$ and a nonzero vector in $V$.

\emph{Cubes.} We fix a grid of cubes of side $K^2$ in $\R^n$ and call its elements $K^2$-cubes. For $\rho\ge K^2$, a \emph{$\rho$-cube} is a cube of a fixed grid whose side lies between $\rho$ and $2\rho$; we choose the grids of $r^{1/2+\eps_0}$-cubes and $r^{1-\alpha}$-cubes that occur below so that every $\rho$-cube is a union of $K^2$-cubes and every $r^{1-\alpha}$-cube is a union of $r^{1/2+\eps_0}$-cubes.

\begin{definition}[Broad norms]\label{def:broad}
Let $f=(f_\tau)_\tau$ be a family, $A\ge0$ an integer and $1\le p<\infty$. For a $K^2$-cube $b$ define
\begin{equation}\label{eq:mu-b}
 \cB_{k,A,p}(f;b)
 =\min_{V_1,\ldots,V_A}
 \max_{\substack{\tau:\,\ang(G(2\tau),V_a)>K^{-1}\\
                         \text{for all }1\le a\le A}}
 \int_b|Ef_\tau|^p,
\end{equation}
where the minimum is taken over all $A$-tuples of $(k-1)$-dimensional linear subspaces of $\R^n$, and the maximum of the empty set is zero. For a union $U$ of $K^2$-cubes, set
\[
 \|Ef\|_{\BL^p_{k,A,K}(U)}^p
 =\sum_{b\subset U}\cB_{k,A,p}(f;b).
\]
\end{definition}

The minimum in \eqref{eq:mu-b} exists because the maximum depends only on the set of caps that satisfy the angle condition, and there are finitely many such sets. Since the subspaces $V_a$ need not be distinct, $\cB_{k,A,p}(f;b)$ is nonincreasing in $A$. For $A=0$ the maximum is taken over all caps.

The broad norm satisfies a triangle inequality; compare \cite[Lemma~4.1]{Guth}.

\begin{lemma}[Broad triangle inequality]\label{lem:triangle}
Let $p\ge1$, let $f=\sum_{j=1}^Jf_j$ be a sum of families, and let $A_1,\ldots,A_J$
be nonnegative integers with $A_1+\cdots+A_J\le A$. Then, for every union $U$ of $K^2$-cubes,
\[
\|Ef\|_{\BL^p_{k,A,K}(U)}
\le
\sum_{j=1}^J
\|Ef_j\|_{\BL^p_{k,A_j,K}(U)}.
\]
\end{lemma}
\begin{proof}
Fix a $K^2$-cube $b\subset U$, and for each $j$ let $V^{(j)}_1,\ldots,V^{(j)}_{A_j}$ be subspaces that attain the minimum in the definition of $\cB_{k,A_j,p}(f_j;b)$. The collection of all these subspaces has at most $A$ elements. If a cap $\tau$ satisfies $\ang(G(2\tau),V^{(j)}_a)>K^{-1}$ for all $j$ and $a$, then by Minkowski's inequality
\[
 \Bigl(\int_b|Ef_\tau|^p\Bigr)^{1/p}\le\sum_{j=1}^J\Bigl(\int_b|Ef_{j,\tau}|^p\Bigr)^{1/p}
 \le\sum_{j=1}^J\cB_{k,A_j,p}(f_j;b)^{1/p}.
\]
Hence $\cB_{k,A,p}(f;b)^{1/p}\le\sum_j\cB_{k,A_j,p}(f_j;b)^{1/p}$. The lemma follows from Minkowski's inequality in $\ell^p$ over the cubes $b\subset U$.
\end{proof}

We shall use Lemma~\ref{lem:triangle} in the following form. If $f=\sum_{j=1}^Lf_j$, then
\begin{equation}\label{eq:triangle-max}
 \|Ef\|_{\BL^p_{k,A,K}(U)}^p\le L^p\max_{1\le j\le L}\|Ef_j\|_{\BL^p_{k,\lfloor A/L\rfloor,K}(U)}^p .
\end{equation}
Thus a decomposition into $L$ pieces costs a factor $L^p$ in the $p$th power of the norm and a factor $L$ in the parameter $A$.

Next we compare the broad norm with the contributions of transverse caps. A $k$-tuple of caps $(\tau_1,\ldots,\tau_k)$ is called \emph{transverse} if
\begin{equation}\label{eq:transverse}
\begin{aligned}
 &\operatorname{dist}(2\tau_i,2\tau_j)\ge\tfrac18K^{-1}\quad\text{for }i\ne j,\qquad\text{and}\\
 &|G(\xi_1)\wedge\cdots\wedge G(\xi_k)|\ge c'_{n,k}K^{-(k-1)}
 \quad\text{for all }\xi_i\in3\tau_i\ \ (1\le i\le k),
\end{aligned}
\end{equation}
where $c'_{n,k}>0$ is the constant determined in the proof of the next lemma. The second condition is imposed on the larger caps $3\tau_i$ because the fine caps that carry the wave packets of $f_{\tau_i}$ may stick out of $2\tau_i$ slightly. For a family $f$ and a fixed tuple, we write $f_i=f_{\tau_i}$.

\begin{lemma}[Transverse caps]\label{eq:extract-bound}
There are constants $c_{n,k},c'_{n,k}>0$ such that the following holds when the caps have radius $c_{n,k}K^{-1}$. Let $A\ge1$, $p\ge1$, let $f$ be a family, and let $X$ be a union of $K^2$-cubes. Then
\[
\begin{aligned}
 \|Ef\|_{\BL^p_{k,A,K}(X)}^p
 &\le
 \sum_{\substack{\tau_1,\ldots,\tau_k\\\mathrm{transverse}}}
 \sum_{b\subset X}\min_{1\le i\le k}\int_b|Ef_i|^p\le
 \sum_{\substack{\tau_1,\ldots,\tau_k\\\mathrm{transverse}}}
 \left(\prod_{i=1}^k\int_X|Ef_i|^p\right)^{1/k}.
\end{aligned}
\]
\end{lemma}
\begin{proof}
The second inequality follows from $\min_ia_i\le\prod_ia_i^{1/k}$ and H\"older's inequality for the sum over $b$. For the first, it suffices to show that for every $K^2$-cube $b$ with $\cB:=\cB_{k,A,p}(f;b)>0$ there is a transverse tuple with $\int_b|Ef_i|^p\ge\cB$ for all $i$. Let $\mathcal S$ be the set of caps $\tau$ with $\int_b|Ef_\tau|^p\ge\cB$, and put $\rho=K^{-1}/\pi$.

\emph{Step 1.} For every affine subspace $P\subset\R^{n-1}$ of dimension $k-2$ there is a cap $\tau\in\mathcal S$ with $\operatorname{dist}(\xi,P)>\rho$ for all $\xi\in2\tau$. Indeed, the vectors $(-2\xi,1)$, $\xi\in P$, span a linear subspace $V_P\subset\R^n$ of dimension $k-1$. Taking $V_1=\cdots=V_A=V_P$ in \eqref{eq:mu-b}, we find a cap $\tau$ with $\ang(G(2\tau),V_P)>K^{-1}$ and $\int_b|Ef_\tau|^p\ge\cB$. Let $\xi\in2\tau$ and let $\eta$ be the point of $P$ nearest to $\xi$. The vectors $(-2\xi,1)$ and $(-2\eta,1)\in V_P$ differ by $(-2(\xi-\eta),0)$, and $|(-2\xi,1)|\ge1$, so the sine of the angle $\phi$ between them is at most $2|\xi-\eta|$. We have $\phi>K^{-1}$. If $\phi\le\pi/2$, then $\phi\le\frac\pi2\sin\phi\le\pi|\xi-\eta|$ and hence $|\xi-\eta|>K^{-1}/\pi=\rho$; if $\phi>\pi/2$, then $2|\xi-\eta|\ge|(-2\xi,1)|\ge1$, and again $|\xi-\eta|>\rho$.

\emph{Step 2.} For points $\xi_1,\ldots,\xi_k\in\R^{n-1}$ write $\operatorname{Vol}(\xi_1,\ldots,\xi_k)=|(\xi_2-\xi_1)\wedge\cdots\wedge(\xi_k-\xi_1)|$. Choose caps $\tau_1,\ldots,\tau_k\in\mathcal S$ whose centers $\xi_i=\xi_{\tau_i}$ maximize $\operatorname{Vol}(\xi_1,\ldots,\xi_k)$. By Step~1 the centers of the caps in $\mathcal S$ do not lie in a $(k-2)$-dimensional affine subspace, so the maximum is positive. Fix $j$ and let $P_j$ be the affine span of the points $\xi_i$, $i\ne j$. Then $\operatorname{Vol}(\xi_1,\ldots,\xi_k)$ is the product of $\operatorname{dist}(\xi_j,P_j)$ and a factor that does not depend on $\xi_j$. By Step~1 there is $\tau\in\mathcal S$ with $\operatorname{dist}(\xi_\tau,P_j)>\rho$, and replacing $\tau_j$ by $\tau$ does not increase the volume; hence
\begin{equation}\label{eq:heights}
 \operatorname{dist}(\xi_j,P_j)>\rho\qquad(1\le j\le k).
\end{equation}
In particular $|\xi_i-\xi_j|>\rho$ for $i\ne j$, and therefore $\operatorname{dist}(2\tau_i,2\tau_j)\ge\rho-4c_{n,k}K^{-1}\ge\frac18K^{-1}$ if $c_{n,k}\le\frac1{24\pi}$.

\emph{Step 3.} Let $\xi_i'\in3\tau_i$, so that $|\xi_i'-\xi_i|\le3c_{n,k}K^{-1}$, and put $w_j=\xi_j-\xi_1$ and $w_j'=\xi_j'-\xi_1'$ for $2\le j\le k$. For $c\in\R^{k-1}$ and each $j$, \eqref{eq:heights} gives $|\sum_ic_iw_i|\ge|c_j|\operatorname{dist}(w_j,\operatorname{span}\{w_i:i\ne j\})\ge|c_j|\rho$, so $|\sum_ic_iw_i|\ge\rho|c|/\sqrt{k-1}$. Since $|w_j'-w_j|\le6c_{n,k}K^{-1}$, we have $|\sum_ic_i(w_i'-w_i)|\le6c_{n,k}K^{-1}\sqrt{k-1}\,|c|$. If $c_{n,k}\le(12\pi k)^{-1}$, it follows that $|\sum_ic_iw_i'|\ge\rho|c|/(2\sqrt{k-1})$ for all $c$, and therefore
\[
 \operatorname{Vol}(\xi_1',\ldots,\xi_k')=|w_2'\wedge\cdots\wedge w_k'|\ge\Bigl(\frac{\rho}{2\sqrt{k-1}}\Bigr)^{k-1}.
\]
Finally, write $g_i=(-2\xi_i',1)$. Then $g_1\wedge\cdots\wedge g_k=g_1\wedge(g_2-g_1)\wedge\cdots\wedge(g_k-g_1)$, the vectors $g_j-g_1=(-2w_j',0)$ lie in the hyperplane $x_n=0$, and $g_1$ has distance at least $1$ from this hyperplane; hence $|g_1\wedge\cdots\wedge g_k|\ge2^{k-1}\operatorname{Vol}(\xi_1',\ldots,\xi_k')$. Since $|g_i|\le\sqrt{17}$, we conclude that $|G(\xi_1')\wedge\cdots\wedge G(\xi_k')|\ge c'_{n,k}K^{-(k-1)}$ with $c'_{n,k}=17^{-k/2}(\pi\sqrt{k-1})^{-(k-1)}$. Thus $(\tau_1,\ldots,\tau_k)$ is transverse, and $\int_b|Ef_{\tau_i}|^p\ge\cB$ because $\tau_i\in\mathcal S$.
\end{proof}

The number of transverse tuples is at most $CK^{k(n-1)}$, so the sums over tuples in Lemma~\ref{eq:extract-bound} may be replaced by maxima at the cost of a factor $K^{O(1)}$.

\subsection{Wave packets and local averages of the \texorpdfstring{$L^2$}{L2}-norm}
\label{sec:packets}
Fix a small parameter $\eps_0>0$ and let $r\ge1$. We use the wave packet decomposition of \cite[Section~3]{Guth} at scale $r$. Let $\Theta_r$ be a finitely overlapping cover of $B^{n-1}(0,2)$ by caps $\theta$ of radius $r^{-1/2}$, with centers $\xi_\theta$; we call them \emph{fine caps}. For $\theta\in\Theta_r$ and $v\in r^{(1+\eps_0)/2}\mathbb Z^{n-1}$ let
\[
 T_{\theta,v}=\bigl\{(x',x_n)\in\R^n:\ |x_n|\le2r,\ |x'+2x_n\xi_\theta+v|\le r^{1/2+\eps_0}\bigr\}.
\]
This is a tube of length $\sim r$ and radius $\sim r^{1/2+\eps_0}$ in the direction $G(\xi_\theta)$. Every $g\in L^2(\R^{n-1})$ with support in $B^{n-1}(0,3/2)$ has a decomposition
\[
 g=\sum_{\theta,v}g_{\theta,v}+g_{\mathrm{err}},\qquad \|g_{\mathrm{err}}\|_\infty\le\RD(r)\|g\|_2,
\]
with the following properties.
\begin{itemize}
\item[(W1)] $\supp g_{\theta,v}\subset\theta$, and $g_{\theta,v}=0$ unless $\theta$ meets $\supp g$.
\item[(W2)] $|Eg_{\theta,v}(x)|\le\RD(r)\|g\|_2$ for $x\in B(0,2r)\setminus T_{\theta,v}$.
\item[(W3)] For every set $\mathcal S$ of pairs $(\theta,v)$,
\[
 \Bigl\|\sum_{(\theta,v)\in\mathcal S}g_{\theta,v}\Bigr\|_2^2
 \sim \sum_{(\theta,v)\in\mathcal S}\|g_{\theta,v}\|_2^2,
 \qquad\text{and}\qquad
 \sum_v\|g_{\theta,v}\|_2^2\lesssim\int_\theta|g|^2\quad\text{for each }\theta.
\]
\end{itemize}
Properties (W1), (W2) and the first part of (W3) are stated in \cite[Section~3]{Guth} for functions on $B^{n-1}$, and the construction applies without change to $B^{n-1}(0,3/2)$; the second part of (W3) follows from the construction and Plancherel's theorem. In \cite{Guth} the tubes are truncated to $B(0,r)$ and (W2) is stated there; the same proof gives (W2) on $B(0,2r)$, because the radius $r^{1/2+\eps_0}$ of the tubes exceeds $r^{1/2}$ by a power of $r$. We need the larger ball in order to work with cubes that meet $B_r$ without being contained in it. For a ball $B_r=B(x_0,r)$ that is not centered at the origin, we apply the decomposition to the function $e^{i(x_0'\cdot\xi+x_{0,n}|\xi|^2)}g(\xi)$; this translates the tubes by $x_0$.

\emph{Packets of a family.} Let $f=(f_\tau)_\tau$ be a family and let $B_r$ be a ball of radius $r$. Applying the decomposition to each $f_\tau$, we obtain functions $f_{\tau,\theta,v}$. We write $2B_r$ for the ball with the same center and radius $2r$. A \emph{packet} of $f$ (at scale $r$, in $B_r$) is a triple $T=(\tau,\theta,v)$ such that $\theta$ meets $\supp f_\tau$ and the tube $T_{\theta,v}$, translated by the center of $B_r$, meets $2B_r$. We write $f_T=f_{\tau,\theta,v}$, $\tau(T)=\tau$ and $\theta(T)=\theta$, and we use the same letter $T$ for the tube. For a set $\mathcal S$ of packets, $f_{\mathcal S}$ denotes the family
\[
 (f_{\mathcal S})_\tau=\sum_{T\in\mathcal S:\,\tau(T)=\tau}f_T .
\]
By (W1), the support of $(f_{\mathcal S})_\tau$ is contained in the $2r^{-1/2}$-neighborhood of $\supp f_\tau$. By (W2), if $\mathcal S'\subset\mathcal S$ are sets of packets and $\mathcal S'$ contains all packets of $\mathcal S$ whose tubes meet a set $U\subset2B_r$, then
\begin{equation}\label{eq:tails}
 |E(f_{\mathcal S})_\tau(x)-E(f_{\mathcal S'})_\tau(x)|\le\RD(r)\|f\|_2\qquad\text{for }x\in U\text{ and all }\tau ,
\end{equation}
and the same holds with $f$ in place of $f_{\mathcal S}$ when $\mathcal S'$ contains all packets of $f$ whose tubes meet $U$.
The number of packets of $f$ is at most $r^{C}$, since $K\le r$ in all our applications.

\emph{Local averages.} For a family $f$ and $r\ge1$, define $M_r(f)\ge0$ by
\[
 M_r(f)^2=\max_{\theta}
 \frac1{|\theta|}\sum_\tau\int_\theta|f_\tau|^2,
\]
where the maximum is taken over all caps $\theta\subset\R^{n-1}$ of radius $r^{-1/2}$. For the family of a function $f\in L^2(B^{n-1})$ this is comparable to $\max_\theta\|f\|_{L^2_{\mathrm{avg}}(\theta)}$, and in particular $M_r(f)\lesssim\|f\|_\infty$. Since the supports lie in $B^{n-1}(0,2)$, which is covered by $\lesssim|\theta|^{-1}$ caps of radius $r^{-1/2}$,
\begin{equation}\label{eq:L2-vs-M}
 \|f\|_2\lesssim M_r(f).
\end{equation}
Moreover, for every set $\mathcal S$ of packets of $f$ at scale $r$ and every $1\le t\le r$,
\begin{equation}\label{eq:packet-average}
 \|f_{\mathcal S}\|_2\lesssim\|f\|_2,
 \qquad
 M_t(f_{\mathcal S}) \lesssim M_r(f).
\end{equation}
The first bound follows from (W3). For the second, let $\vartheta$ be a cap of radius $t^{-1/2}\ge r^{-1/2}$. By (W1) only the packets with $\theta(T)\cap\vartheta\ne\varnothing$ contribute to $f_{\mathcal S}$ on $\vartheta$, so by (W3)
\[
 \sum_\tau\int_\vartheta|(f_{\mathcal S})_\tau|^2\lesssim\sum_\tau\ \sum_{\theta:\,\theta\cap\vartheta\ne\varnothing}\int_\theta|f_\tau|^2
 \le\sum_{\theta:\,\theta\cap\vartheta\ne\varnothing}|\theta|M_r(f)^2\lesssim|\vartheta|M_r(f)^2 .
\]

\emph{Uniform sets of packets.} A set $\T$ of packets of $f$ is called \emph{uniform with parameters $(a,m)$} if
\begin{equation}\label{eq:uniform}
\begin{aligned}
 &a\le\|f_T\|_2<2a&&\text{for all }T\in\T,\\
 &\#\{T\in\T:\theta(T)=\theta\}\in\{0\}\cup[m,2m)&&\text{for all }\theta\in\Theta_r .
\end{aligned}
\end{equation}
For such a set, (W3) gives
\begin{equation}\label{eq:am}
\begin{aligned}
 a^2m&\lesssim r^{-(n-1)/2}
 M_r(f)^2, \qquad a^2\#\T\lesssim\|f\|_2^2 .
\end{aligned}
\end{equation}
Indeed, if $\theta$ is a fine cap that contains a packet of $\T$, then
\[
 a^2m\le\sum_{T\in\T:\,\theta(T)=\theta}\|f_T\|_2^2\lesssim\sum_\tau\int_\theta|f_\tau|^2\le|\theta|M_r(f)^2 ,
\]
and the second bound follows in the same way by summing over all fine caps. A subset of a uniform set has at most $2m$ packets in each fine cap, and all its packets satisfy $\|f_T\|_2\sim a$.

\subsection{An \texorpdfstring{$L^2$}{L2} estimate and refined decoupling}\label{subsec:two-analytic}
We now state the two estimates \eqref{eq:intro-L2} in the form in which they will be used. Let $f$ be a family, let $B_r$ be a ball of radius $r$, and let $\T$ be a uniform set of packets of $f$ with parameters $(a,m)$. Let $\tau_i$ be a cap, let $\T_i$ be a set of packets of $\T$ with $\tau(T)=\tau_i$, and put
\[
 f_i=\sum_{T\in\T_i}f_T,
 \qquad N_i=\#\T_i ;
\]
by (W3), $\|f_i\|_2^2\sim a^2N_i$. Let $X\subset2B_r$ be a union of $r^{1/2+\eps_0}$-cubes. For each $T\in\T_i$, let
$Y_i(T)$ be the union of the cubes of $X$ that meet $T$, and define the average density $\lambda_i$ by
\[
 \sum_{T\in\T_i}|Y_i(T)|
 =\lambda_i\sum_{T\in\T_i}|T|.
\]

\emph{The $L^2$ estimate.} We have
\begin{equation}\label{eq:L2}
 \int_X|Ef_i|^2
 \lesssim r^{O(\eps_0)}\lambda_i r\,\|f_i\|_2^2+\RD(r)\|f\|_2^2.
\end{equation}
To see this, let $Q\subset X$ be one of the cubes and let $\T_{i,Q}$ be the set of $T\in\T_i$ that meet $Q$. By \eqref{eq:tails}, $Ef_i$ agrees on $Q$ with $E\bigl(\sum_{T\in\T_{i,Q}}f_T\bigr)$ up to $\RD(r)\|f\|_2$. The elementary bound $\int_{B_\rho}|Eg|^2\lesssim\rho\|g\|_2^2$ \cite[Lemma~3.2]{Guth} and (W3) give $\int_Q|Ef_i|^2\lesssim r^{1/2+\eps_0}a^2\,\#\T_{i,Q}+\RD(r)\|f\|_2^2$. Summing over $Q$ and using
\[
 \sum_{Q\subset X}\#\T_{i,Q}=\sum_{T\in\T_i}\frac{|Y_i(T)|}{|Q|}=\lambda_iN_i\frac{|T|}{|Q|}\lesssim\lambda_iN_i\,r^{1/2},
\]
we obtain \eqref{eq:L2}.

\emph{Refined decoupling.} Set $q=2(n+1)/(n-1)$, and suppose that every cube of $X$ meets at most $\mu_i$
tubes of $\T_i$. Then, for every $\eps>0$,
\begin{equation}\label{eq:RD}
\begin{aligned}
 \int_X|Ef_i|^q
 &\lesssim_{\eps}
 r^{\eps+O(\eps_0)}\mu_i^{2/(n-1)}a^qN_i+\RD(r)\|f\|_2^q\\
 &\lesssim
 r^{-1+\eps+O(\eps_0)}
 (\mu_i/m)^{2/(n-1)}\|f_i\|_2^2
 M_r(f)^{q-2}+\RD(r)\|f\|_2^q.
\end{aligned}
\end{equation}
The first inequality follows from the refined decoupling theorem \cite[Theorem~4.2]{GIOW} (see also \cite[Theorem~0.3]{WW}): if the norms $\|Ef_T\|_{L^q(w_{B_r})}$, $T\in\T_i$, are comparable to each other, and every cube of side $r^{1/2}$ contained in $X$ meets at most $M$ tubes of $\T_i$, then
\[
 \int_X|Ef_i|^q\lesssim_\eps r^{\eps}M^{2/(n-1)}\sum_{T\in\T_i}\|Ef_T\|_{L^q(w_{B_r})}^q ,
\]
where $w_{B_r}$ is a weight that is comparable to $1$ on $B_r$ and decays rapidly away from $B_r$. The theorem is stated for sets $X$ contained in a ball of radius $r$; we apply it on each of $O(1)$ balls of radius $r$ that cover $2B_r$. In our situation we may take $M=\mu_i$. Moreover $\|Ef_T\|_\infty\le\|f_T\|_1\lesssim r^{-(n-1)/4}a$ by (W1) and $\|Ef_T\|_{L^2(w_{B_r})}^2\lesssim ra^2$, so that $\|Ef_T\|_{L^q(w_{B_r})}^q\lesssim a^q$, because $(q-2)(n-1)/4=1$. To arrange that the norms are comparable, we split $\T_i$ into $O(\log r)$ subsets according to the dyadic value of $\|Ef_T\|_{L^q(w_{B_r})}$, discarding the packets for which this norm is less than $r^{-C}a$, and apply the theorem to each subset; this costs a factor $(\log r)^{q}$, which is absorbed in $r^\eps$. The second inequality in \eqref{eq:RD} follows from $a^qN_i\sim a^{q-2}\|f_i\|_2^2$ and from the first bound in \eqref{eq:am}, since $(q-2)/2=2/(n-1)$.

\emph{Rescaling to unit length.} The geometric estimates of Section~\ref{subsec:geom-tools} will be applied after the change of variables $x\mapsto x/r$ (preceded by a translation that moves the center of $B_r$ to the origin). Let $s$ be the side length of the $r^{1/2+\eps_0}$-cubes and put
\[
 \delta_*=s/r,\qquad\text{so that}\qquad r^{-1/2+\eps_0}\le\delta_*\le2r^{-1/2+\eps_0}.
\]
The rescaled tubes $T_{\theta,v}$ are $\delta_*$-tubes in the sense of Definition~\ref{def:tubes}, with direction $G(\xi_\theta)$, and the $r^{1/2+\eps_0}$-cubes become the cubes of the grid $\cG_{\delta_*}$ (after a translation of the grid); thus $T\mapsto Y_i(T)$ becomes a shading in the sense of Definition~\ref{def:shading}. Packets with the same fine cap $\theta$ have parallel tubes, and the directions $G(\xi_\theta)$ that belong to one direction class at scale $\delta_*$ come from at most $Cr^{(n-1)\eps_0}$ fine caps. Consequently, if $\T$ is a uniform set of packets with parameters $(a,m)$, then the rescaled tubes of every subset of $\T$ form an
\begin{equation}\label{eq:mstar}
 m_*\text{-parallel family, with}\quad m_*\le Cr^{(n-1)\eps_0}m .
\end{equation}

\section{A multilinear incidence estimate}\label{sec:newframework}
Fix $n,k$ as in \eqref{eq:k-range} and set $d=d_{n,k}$. In this section we deduce from Lemma~\ref{thm:zahl-input} a bound for the multiplicity of $k$ transverse families of tubes. It has the form \eqref{eq:intro-ab} with $(a,b)=(n-d,d-1)$. Two features of the statement are dictated by the application: the families are $m$-parallel rather than direction-separated, and we do not assume that they have comparable cardinalities. The estimate will be combined with the hairbrush bound and refined decoupling in Section~\ref{sec:proof}.

\begin{lemma}\label{lem:geometry}
Let $\eps>0$, $\nu>0$, and let $m\ge1$ be an integer. Let $\T_1,\ldots,\T_k$ be families of $\delta$-tubes in $\R^n$, let $Y_i$ be a shading of $\T_i$ for each $i$, and let $E\subset B(0,C_n)$ be a union of cubes of $\cG_\delta$. Assume the following.
\begin{itemize}
\item[(a)] The union $\T=\T_1\cup\cdots\cup\T_k$ is $m$-parallel.
\item[(b)] We have $|v(T)-v(T')|\ge\delta$ whenever $T$ and $T'$ belong to different families, and $|v(T_1)\wedge\cdots\wedge v(T_k)|\ge\nu$ for every choice of $T_i\in\T_i$.
\item[(c)] There are numbers $\mu_i>0$ such that $\mu_i\le\mu_{Y_i}(x)\le2\mu_i$ for all $x\in E$ and all $1\le i\le k$.
\end{itemize}
Define the average density $\lambda_i$ of $Y_i$ on $E$ by $\sum_{T\in\T_i}|Y_i(T)\cap E|=\lambda_i\sum_{T\in\T_i}|T|$, and set
\[
 \begin{gathered}
 N_i=\#\T_i,\qquad N=\sum_iN_i,\qquad
 Z=\frac{\delta^{n-1}N}{m},\\
 \zeta=\frac{(N_1\cdots N_k)^{1/k}}N,\qquad
 \Lambda=(\lambda_1\cdots\lambda_k)^{1/k},\qquad
 \mu=(\mu_1\cdots\mu_k)^{1/k}.
 \end{gathered}
\]
If $\min_i \mu_i/m\ge C_{n,k}\log(2/\delta)$ for a sufficiently large constant $C_{n,k}$, then
\begin{equation}\label{eq:geom-conclusion}
 \frac\mu m\lesssim_{\eps}
 \nu^{-1}\delta^{-(n-d)-\eps}
 \Lambda^{-(d-1)}\zeta^{-(d-1)}Z^{(d-1)/(n-1)}.
\end{equation}
\end{lemma}

In \eqref{eq:geom-conclusion}, the factor $\nu^{-1}$ is the cost of transversality, and $\zeta^{-(d-1)}$ is a loss that occurs when the families have very different cardinalities; note that $\zeta\le1/k$ by the inequality between the arithmetic and geometric means. In the application the loss in $\zeta$ will be compensated by a gain in the $L^2$ norms; see the remark after \eqref{multilinear-norm}. By \eqref{eq:card-parallel} we have $Z\lesssim1$, so the last factor is a gain, which we shall not use. The lower bound on $\mu_i/m$ is needed for the random sampling in the proof.

\begin{proof}
\emph{Step 1: sampling.} For $\omega\in\Omega_\delta$, the direction class $\T(\omega)$ of $\omega$ in $\T$ has at most $m$ elements. Independently for each $\omega$, we select at most one tube of $\T(\omega)$ at random, in such a way that each tube of $\T(\omega)$ is selected with probability $1/m$. Let $\T'$ be the set of selected tubes and $\T'_i=\T'\cap\T_i$. The family $\T'$ is $1$-parallel.

Fix a cube $Q$ of $\cG_\delta$ contained in $E$ and an index $i$. The number
\[
 \mu_i'(Q)=\#\{T\in\T_i':Q\subset Y_i(T)\}
\]
is a sum of independent random variables with values in $\{0,1\}$, one for each direction class, and its expectation is $\mu_{Y_i}(Q)/m\ge\mu_i/m\ge C_{n,k}\log(2/\delta)$. By Chernoff's inequality, the probability that $\mu_i'(Q)<\mu_i/(2m)$ is at most $\exp(-\mu_i/(8m))\le(\delta/2)^{C_{n,k}/8}$. There are $\lesssim\delta^{-n}$ cubes of $\cG_\delta$ in $E$ and $k$ indices $i$. Hence, if $C_{n,k}$ is large enough, then with probability greater than $3/4$,
\begin{equation}\label{eq:sampled-mult}
 \mu_i'(Q)\ge\frac{\mu_i}{2m}\qquad\text{for all cubes }Q\subset E\text{ and all }1\le i\le k.
\end{equation}
The expected number of selected tubes is $N/m$, so $\#\T'\le8N/m$ with probability at least $7/8$, by Markov's inequality. We fix a selection $\T'$ for which \eqref{eq:sampled-mult} holds and $\#\T'\le8N/m$.

\emph{Step 2: Zahl's theorem.} We use the partition \eqref{eq:coloring} of $\Omega_\delta$. For $\mathbf j=(j_1,\ldots,j_k)\in\{1,\ldots,c_n\}^k$, let
\[
 \T'_{\mathbf j}=\bigcup_{i=1}^k\bigl\{T\in\T_i':\omega(T)\in\Omega_\delta^{(j_i)}\bigr\}.
\]
Each family $\T'_{\mathbf j}$ is direction-separated: two of its tubes that belong to the same $\T_i'$ have directions at distance at least $\delta$ by the argument that precedes \eqref{eq:one-parallel}, and two of its tubes that belong to different families have directions at distance at least $\delta$ by hypothesis~(b). For a tube $T$ let $\widetilde T=N_{\sqrt n\,\delta}(T)$; this is again a $\delta$-tube, with a larger constant in Definition~\ref{def:tubes}. If $Q\subset Y_i(T)$, then $Q$ meets $T$ and therefore $Q\subset\widetilde T$. Let
\[
 F_{\mathbf j}=\sum_{T_1,\ldots,T_k\in\T'_{\mathbf j}}\1_{\widetilde T_1}\cdots\1_{\widetilde T_k}\,|v(T_1)\wedge\cdots\wedge v(T_k)|^{k/d}.
\]
Every $k$-tuple $(T_1,\ldots,T_k)\in\T_1'\times\cdots\times\T_k'$ consists of tubes of $\T'_{\mathbf j}$ for exactly one $\mathbf j$. Hence, by \eqref{eq:sampled-mult} and hypothesis~(b), for every $x\in E$,
\[
\begin{aligned}
 \nu^{k/d}\Bigl(\frac{\mu}{2m}\Bigr)^k&=\nu^{k/d}\prod_{i=1}^k\frac{\mu_i}{2m}\\
 &\le\sum_{T_1\in\T_1',\ldots,T_k\in\T_k'}\1_{\widetilde T_1}(x)\cdots\1_{\widetilde T_k}(x)\,|v(T_1)\wedge\cdots\wedge v(T_k)|^{k/d}
 \le\sum_{\mathbf j}F_{\mathbf j}(x).
\end{aligned}
\]
Put $s=d/(k(d-1))$. We take the $L^s(E)$ quasi-norm of both sides and apply Lemma~\ref{thm:zahl-input} to each of the $c_n^k$ families $\T'_{\mathbf j}$, with the tubes $\widetilde T$. Since $\delta^{n-1}\#\T'_{\mathbf j}\le8\delta^{n-1}N/m=8Z$, this gives
\[
 \nu^{k/d}\Bigl(\frac{\mu}{2m}\Bigr)^k|E|^{1/s}\lesssim_\eps\Bigl(\delta^{1-n/d-\eps}Z^{\frac{n(d-1)}{(n-1)d}}\Bigr)^k .
\]
Taking $k$th roots and raising the result to the power $d/(d-1)$, we obtain
\begin{equation}\label{eq:geom-mid}
 \nu^{1/(d-1)}(\mu/m)^{d/(d-1)}|E|
 \lesssim_{\eps}\delta^{-(n-d)/(d-1)-\eps}
 Z^{n/(n-1)},
\end{equation}
with a different value of $\eps$.

\emph{Step 3: conclusion.} By hypothesis~(c) and the definition of $\lambda_i$, we have $\mu_i|E|\sim\lambda_i\delta^{n-1}N_i$ for each $i$. Taking geometric means gives
\[
 \mu|E|\sim\Lambda\zeta\delta^{n-1}N
   =\Lambda\zeta m Z.
\]
We substitute $|E|\sim\Lambda\zeta mZ/\mu$ into \eqref{eq:geom-mid}, divide by $\Lambda\zeta Z$, and raise the result to the power $d-1$. This proves \eqref{eq:geom-conclusion}.
\end{proof}

\section{Proof of the main theorem}\label{sec:proof}
\subsection{The broad two-ends induction}\label{sec:induction}
We prove the broad estimate by a two-ends induction for the broad norms of Definition~\ref{def:broad}, combining Lemma~\ref{lem:geometry} with the hairbrush bound and refined decoupling in the two-ends case. The outer scale $R$, and with it $K$, remains fixed while the intermediate radius $r$ decreases. The broad parameter $A$ is allowed to be a small power of $r$; this pays for the dyadic decompositions, each of which divides $A$ by a power of $\log r$. Compare \cite[Proposition~4.13]{DW}.

\begin{proposition}\label{prop:broad}
Fix $n,k$ satisfying \eqref{eq:k-range} and write $p_b=p_b(n,k)$. For every $\eta>0$ there is $\kappa_0=\kappa_0(n,k,\eta)>0$ such that the following holds for all $0<\kappa\le\kappa_0$ and all $a_0>0$. Let $R\ge1$, let $K$ be a dyadic number with $R^\kappa\le K<2R^\kappa$, and let $A\ge R^{a_0}$ be an integer. Then, for every $f\in L^2(B^{n-1})$,
\begin{equation}\label{eq:broad-main}
 \|E f\|_{\BL^{p_b}_{k,A,K}(B_R)}^{p_b}
 \lesssim_\eta R^{2\eta}\| f\|_2^2M_R(f)^{p_b-2}.
\end{equation}
The implicit constant depends on $n,k,\eta,\kappa$ and $a_0$, but not on $R$, $A$ or $f$.
\end{proposition}

\begin{proof}
Throughout the proof, $p=p_b$ and $q=2(n+1)/(n-1)$.

\step{Parameter choices.}
Choose $\sigma>0$ with $n\sigma<\eta/4$, then $\alpha\in(0,\frac12)$ so small
that the factors $r^{O(\alpha)}$ below are at most $r^{\eta/8}$.
Choose $b_0\ll\alpha\eta$ and $\eps_0\ll b_0$, and take the losses $\eps$ in Lemma~\ref{thm:GWX-input}, in Lemma~\ref{lem:geometry} and in \eqref{eq:RD} so small that the factors $r^{\eps+O(\eps_0)}$ below are at most $r^{\eta/8}$. Choose $\kappa_0$ so that $4\kappa_0<\sigma$ and so that the factors $K^{O(1)}$ below are at most $R^{\eta/4}$ when $\kappa\le\kappa_0$. Since $K<2R^\kappa$, the condition $4\kappa<\sigma$ ensures that $K^2\le r^{1/2}$ whenever $r\ge R^{\sigma}$ and $R$ is large, so that the $r^{1/2+\eps_0}$-cubes are unions of $K^2$-cubes; it also gives $r^{-1/2}\le c_{n,k}K^{-1}$ and $K\le r$ for such $r$. All parameters are fixed independently of $R$.

\step{The inductive statement.}
We shall prove that there is a constant $C_*$ such that, for all $R$ and $K$ as in the proposition, all $1\le r\le R$, all integers $A\ge r^{a_0}$, all balls $B_r$ of radius $r$, and all families $f=(f_\tau)_\tau$,
\begin{equation}\label{eq:IH}
 \|Ef\|_{\BL^{p_b}_{k,A,K}(B_r)}^{p_b}
 \le C_* R^\eta r^\eta \|f\|_2^2M_r(f)^{p_b-2}.
\end{equation}
For $r=R$ and the family of a function $f\in L^2(B^{n-1})$ this is \eqref{eq:broad-main}. We may assume that $R$ is larger than a constant $R_0$ that depends on the parameters, because for $R\le R_0$ the argument for the base case below gives \eqref{eq:IH} with a constant depending on $R_0$.

The proof is by induction on $r$: we treat $r\le R^\sigma$ directly, and for $r>R^\sigma$ we deduce \eqref{eq:IH} at radius $r$ from \eqref{eq:IH} at radius $r_1=r^{1-\alpha}$. The scales $R$ and $K$ are not changed. Starting from any $r\le R$, the sequence $r,r^{1-\alpha},r^{(1-\alpha)^2},\ldots$ falls below $R^\sigma$ after at most $J$ steps, where $J$ depends only on $\sigma$ and $\alpha$.

In the inductive step, \eqref{eq:IH} at radius $r_1$ is applied to families of the form $f_{\mathcal S}$, where $\mathcal S$ is a set of packets of $f$ at scale $r$. By (W1) their components are supported in the $2r^{-1/2}$-neighborhoods of the supports of the $f_\tau$. Since at most $J$ steps occur and packets are formed only at scales $r>R^\sigma$, the supports grow by at most $2JR^{-\sigma/2}$ in total; this is less than $c_{n,k}K^{-1}$ because $4\kappa<\sigma$. Hence, if we start from a function $f\in L^2(B^{n-1})$, for which $\supp f_\tau\subset\tau$, then all families that occur have their components supported in $2\tau$, as required in Section~\ref{sec:broad}. We do not record this enlargement in the notation.

\step{Base case.}
Let $r\le R^\sigma$. We have $|Ef_\tau|\le\|f_\tau\|_1\lesssim\|f\|_2$ pointwise, so each $K^2$-cube $b$ and each cap $\tau$ satisfy $\int_b|Ef_\tau|^{p}\lesssim K^{2n}\|f\|_2^{p}$. Summing over the cubes $b\subset B_r$, we get
\[
 \|Ef\|_{\BL^{p}_{k,A,K}(B_r)}^{p}\lesssim r^n\|f\|_2^{p}.
\]
Since $\|f\|_2\lesssim M_r(f)$ by \eqref{eq:L2-vs-M} and $r^n\le R^{n\sigma}\le R^{\eta/4}$, this is bounded by the right-hand side of \eqref{eq:IH} if $C_*$ is large enough.

\step{Negligible terms.}
From now on let $R^\sigma<r\le R$, and assume \eqref{eq:IH} at radius $r_1=r^{1-\alpha}$. Fix a ball $B_r$, an integer $A\ge r^{a_0}$ and a family $f$. Let $C_1$ be a large constant. We may assume that
\begin{equation}\label{eq:nontrivial}
 \|Ef\|_{\BL^{p}_{k,A,K}(B_r)}^{p}\ge r^{-C_1}\|f\|_2^{p},
\end{equation}
since otherwise \eqref{eq:IH} follows from \eqref{eq:L2-vs-M}. We call a quantity \emph{negligible} if it is at most $r^{-3C_1}\|f\|_2^{p}$. If a family $h$ satisfies $|Eh_\tau|\le r^{-C_2}\|f\|_2$ on $B_r$ for all $\tau$, with $C_2$ sufficiently large in terms of $C_1$, then $\|Eh\|_{\BL^p_{k,0,K}(B_r)}^p$ is negligible; by Lemma~\ref{lem:triangle}, such a family may be removed from $f$ without changing the broad parameter. This applies in particular to the errors in \eqref{eq:tails}. Negligible terms remain much smaller than the left-hand side of \eqref{eq:nontrivial} after multiplication by the factors $K^{O(1)}r^{O(\alpha)}(\log r)^{O(1)}$ that occur below, and we shall omit them.

\step{Step 1: uniform packets.}
Decompose $f$ into packets at scale $r$ in $B_r$, as in Section~\ref{sec:packets}. There are at most $r^{C}$ packets, so the packets with $\|f_T\|_2<r^{-C_2-C}\|f\|_2$ together form a family of the kind just described, and we discard them. We sort the remaining packets into classes: for dyadic numbers $a$ and $m$, the packet $T$ belongs to the class $(a,m)$ if $a\le\|f_T\|_2<2a$ and if the number of packets $T'$ with $a\le\|f_{T'}\|_2<2a$ and $\theta(T')=\theta(T)$ lies in $[m,2m)$. Each class is a uniform set of packets with parameters $(a,m)$ in the sense of \eqref{eq:uniform}. Since $r^{-C}\|f\|_2\lesssim a\lesssim\|f\|_2$ and $1\le m\le r^{C}$, the number $L_1$ of classes is $\lesssim(\log r)^2$. By \eqref{eq:triangle-max} there is a class, which we denote by $\T^0$ and whose parameters we denote by $(a,m)$, such that
\begin{equation}\label{eq:step1}
 \|Ef\|_{\BL^{p}_{k,A,K}(B_r)}^{p}\lesssim(\log r)^{2p}\,\|Ef_{\T^0}\|_{\BL^{p}_{k,A_1,K}(B_r)}^{p},
 \qquad A_1=\lfloor A/L_1\rfloor .
\end{equation}

\step{Step 2: uniform cubes.}
Consider the $r^{1/2+\eps_0}$-cubes $Q$ that meet $B_r$; they are contained in $2B_r$. For such a cube write
\[
 \mathcal B(Q)
 :=\|Ef_{\T^0}\|_{\BL^{p_b}_{k,A_1,K}(Q\cap B_r)}^{p_b},
\]
where $Q\cap B_r$ stands for the union of the $K^2$-cubes $b\subset Q$ that are contained in $B_r$; thus $\sum_Q\mathcal B(Q)=\|Ef_{\T^0}\|_{\BL^{p}_{k,A_1,K}(B_r)}^{p}$. As in the base case, $\mathcal B(Q)\lesssim r^{n}\|f\|_2^{p}$ for every $Q$. The cubes with $\mathcal B(Q)<r^{-4C_1}\|f\|_2^{p}$ contribute a negligible amount in total. By dyadic pigeonholing, there are a collection $\mathcal Q$ of cubes and a number $H\ge r^{-4C_1}\|f\|_2^{p}$ such that, with
\[
 X=\bigcup_{Q\in\mathcal Q}Q,
\]
we have
\begin{equation}\label{eq:step2}
 H\le\mathcal B(Q)<2H
 \quad(Q\in\mathcal Q),
 \qquad
 \sum_{Q\in\mathcal Q}\mathcal B(Q)
 \gtrsim(\log r)^{-1}
 \|Ef_{\T^0}\|_{\BL^{p_b}_{k,A_1,K}(B_r)}^{p_b}.
\end{equation}

\step{Step 3: the density $\lambda$ and the number $\beta$ of slabs.}
Partition the $x_n$-axis into intervals $J$ whose length is the side length of the $r_1$-cubes, aligned with the grid, and call $\R^{n-1}\times J$ the \emph{slab} of $J$. Each cube $Q$ lies in exactly one slab, and we write $J(Q)$ for the corresponding interval. The number of slabs that meet $B_r$ is $\lesssim r^\alpha$.
For each $T\in\T^0$ and each interval $J$, define
\[
 Y_J(T)
 :=\bigcup_{\substack{Q\in\mathcal Q:\,J(Q)=J\\Q\cap T\ne\varnothing}}Q .
\]
For dyadic numbers $\lambda,\beta>0$,
set
\[
 \mathcal J_\lambda(T)
 :=\{J:\lambda|T|\le |Y_J(T)|<2\lambda|T|\},
 \qquad
 \mathbb T_{\lambda,\beta}
 :=\{T\in\T^0:\beta\le\#\mathcal J_\lambda(T)<2\beta\}.
\]
Thus $\mathcal J_\lambda(T)$ is the set of slabs in which $X$ has density about $\lambda$ on $T$, relative to the whole tube, and $\beta$ counts these slabs. If $Y_J(T)$ is nonempty, then $|Q|\le|Y_J(T)|\lesssim r^{-\alpha}|T|$, so only $O(\log r)$ values of $\lambda$ occur; moreover $1\le\beta\lesssim r^\alpha$. Hence the number $L_2$ of pairs $(\lambda,\beta)$ that occur is $\lesssim(\log r)^2$. We say that $T\in\T_{\lambda,\beta}$ is \emph{active} on $J$ if $J\in\mathcal J_\lambda(T)$, and we write
\[
 \T_{\lambda,\beta,J}=\{T\in\T_{\lambda,\beta}:J\in\mathcal J_\lambda(T)\}.
\]

Fix $Q\in\mathcal Q$ and put $J=J(Q)$. If the tube of $T\in\T^0$ meets $Q$, then $Y_J(T)\supset Q$ is nonempty, so there is exactly one $\lambda$ with $J\in\mathcal J_\lambda(T)$, and then exactly one $\beta$ with $T\in\T_{\lambda,\beta}$. In other words, the sets $\T_{\lambda,\beta,J}$, for the various pairs $(\lambda,\beta)$, are pairwise disjoint, and their union contains every packet of $\T^0$ whose tube meets $Q$. By \eqref{eq:tails}, on $Q$ the family $f_{\T^0}$ differs from $\sum_{\lambda,\beta}f_{\T_{\lambda,\beta,J}}$ by a negligible family. Choose
\[
 A'=\left\lfloor\frac{A_1}{L_2}\right\rfloor,
 \qquad\text{and let}\qquad
 \mathcal B_{\lambda,\beta}(Q):=\|Ef_{\T_{\lambda,\beta,J(Q)}}\|_{\BL^{p_b}_{k,A',K}(Q\cap B_r)}^{p_b}.
\]
By \eqref{eq:triangle-max},
\[
 \mathcal B(Q)
 \lesssim(\log r)^{2p}
 \max_{\lambda,\beta}\mathcal B_{\lambda,\beta}(Q).
\]
Hence, for every $Q\in\mathcal Q$, we can choose
$\lambda(Q),\beta(Q)$ such that
\[
 \mathcal B_{\lambda(Q),\beta(Q)}(Q)
 \gtrsim(\log r)^{-2p}H.
\]
Choose a dyadic number $h(Q)$ satisfying
\[
 h(Q)\le
 \mathcal B_{\lambda(Q),\beta(Q)}(Q)<2h(Q).
\]
As in the base case, $\mathcal B_{\lambda,\beta}(Q)\lesssim r^{n}\|f\|_2^{p_b}$; here we use the first bound in \eqref{eq:packet-average}.
Together with the lower bound for $H$, this implies
\[
 (\log r)^{-2p}H\lesssim h(Q)\lesssim r^{n+4C_1} H.
\]
Thus $h(Q)$ takes only $O(\log r)$ dyadic values.

Pigeonholing the cubes once more, choose a fixed triple
$(\lambda,\beta,h_0)$ and a subcollection
$\mathcal Q_2\subset\mathcal Q$ such that
\[
 \#\mathcal Q_2\gtrsim(\log r)^{-3}\#\mathcal Q,
 \qquad
 (\lambda(Q),\beta(Q),h(Q))=(\lambda,\beta,h_0)
 \quad(Q\in\mathcal Q_2).
\]
Define
\[
 X_2:=\bigcup_{Q\in\mathcal Q_2}Q.
\]
Since the cubes have equal volume,
\begin{equation}\label{eq:refinement-size}
 |X_2|\gtrsim(\log r)^{-3}|X|.
\end{equation}
Moreover,
\[
 h_0\le\mathcal B_{\lambda,\beta}(Q)<2h_0,
 \qquad h_0\gtrsim(\log r)^{-2p}H
 \qquad(Q\in\mathcal Q_2).
\]
For a union $U$ of cubes of $\mathcal Q_2$, write
\[
 \mathcal B_{\lambda,\beta}(U)
 :=\sum_{Q\subset U}\mathcal B_{\lambda,\beta}(Q).
\]
Then $\mathcal B_{\lambda,\beta}(X_2)\ge h_0\,\#\cQ_2\gtrsim(\log r)^{-C}H\,\#\cQ\gtrsim(\log r)^{-C}\sum_{Q\in\cQ}\cB(Q)$, and by \eqref{eq:step1} and \eqref{eq:step2},
\begin{equation}\label{eq:step3}
 \|Ef\|_{\BL^{p_b}_{k,A,K}(B_r)}^{p_b}\lesssim(\log r)^{C}\,\mathcal B_{\lambda,\beta}(X_2).
\end{equation}
Thus the passage from $\cQ$ to $\cQ_2$ and from $f$ to the families $f_{\T_{\lambda,\beta,J}}$ costs only powers of $\log r$, both in the volume and in the broad norm, and the quantities $\cB_{\lambda,\beta}(Q)$, $Q\in\cQ_2$, are comparable to each other. The broad parameter has been divided by $L_1L_2\lesssim(\log r)^4$: for sufficiently large $r$,
\begin{equation}\label{eq:A-prime}
 A'\ge\frac{A}{L_1L_2}-2\gtrsim\frac{r^{a_0}}{(\log r)^{4}}
 \ge r^{(1-\alpha)a_0}=r_1^{a_0}.
\end{equation}

\step{Step 4: the shading.}
From now on the pair $(\lambda,\beta)$ is fixed, and we write $\T=\T_{\lambda,\beta}$ and $\T_J=\T_{\lambda,\beta,J}$. For $T\in\mathbb T$, define the shading
\[
 Y(T):=\bigcup_{J\in\mathcal J_\lambda(T)}Y_J(T).
\]
Since $T$ is active on at least $\beta$ and fewer than $2\beta$ intervals $J$, and $\lambda|T|\le|Y_J(T)|<2\lambda|T|$ for each of them,
\begin{equation}\label{eq:lambda0}
 \lambda_0|T|\le|Y(T)|<4\lambda_0|T|,
 \qquad \lambda_0:=\beta\lambda.
\end{equation}
These shadings are formed with the cubes of $\cQ$, not only with those of $\cQ_2$. If $Q\in\cQ_2$, then some tube of $\T_{J(Q)}$ meets $Q$, since otherwise $\cB_{\lambda,\beta}(Q)$ would be negligible by \eqref{eq:tails}; hence $Q\subset Y(T)$ for some $T\in\T$. Therefore
\[
 X_2\subset
 \bigcup_{T\in\mathbb T}Y(T)
 \subset X.
\]
We now redefine
\[
 X:=\bigcup_{T\in\mathbb T}Y(T).
\]
The estimate \eqref{eq:refinement-size} remains valid with this
new $X$. We distinguish two cases, according to the size of $\beta$.

\step{Case 1: The non-two-ends case.}
Suppose $\beta\le r^{b_0}$. Then each packet is active on few slabs, and we apply the induction hypothesis in each slab. Cover $B_r$ by balls $B_j$ of radius $r_1$ in such a way that every $K^2$-cube contained in $B_r$ lies in at least one $B_j$ and every point lies in at most $C$ of them. For a ball $B_j$ and an interval $J$ whose slab meets $B_j$, let $\mathcal S_{j,J}$ be the set of packets $T\in\T_J$ whose tubes meet $B_j$, and let $f_{j,J}=f_{\mathcal S_{j,J}}$. By \eqref{eq:tails}, the family $f_{\T_J}$ agrees with $f_{j,J}$ on $B_j$ up to a negligible family. Every $K^2$-cube of $X_2\cap B_r$ lies in some ball $B_j$ and in the slab of some interval $J$, and this slab meets $B_j$. Hence
\[
 \mathcal B_{\lambda,\beta}(X_2)
 \lesssim
 \sum_{j,J}\|Ef_{j,J}\|_{\BL^{p_b}_{k,A',K}(B_j)}^{p_b}.
\]
A packet $T\in\T$ is active on fewer than $2\beta\le2r^{b_0}$ intervals $J$. For each of them, the balls $B_j$ that meet both the tube of $T$ and the slab of $J$ lie within distance $O(r_1)$ of the part of the tube in that slab, which has diameter $\lesssim r_1$; so there are $O(1)$ such balls. Thus each packet belongs to $O(r^{b_0})$ of the sets $\mathcal S_{j,J}$, and by (W3) and \eqref{eq:packet-average},
\[
\sum_{j,J}\|f_{j,J}\|_2^2
\lesssim r^{b_0}\sum_{T\in\T}\|f_T\|_2^2\lesssim r^{b_0}\|f\|_2^2,
\qquad
M_{r_1}(f_{j,J})\lesssim M_r(f).
\]
As $A'\ge r_1^{a_0}$ by \eqref{eq:A-prime}, the induction hypothesis
\eqref{eq:IH} at radius $r_1$ applies to each $f_{j,J}$ on $B_j$ and yields
\[
\begin{aligned}
\mathcal B_{\lambda,\beta}(X_2)
&\lesssim
C_* R^\eta r_1^\eta
\sum_{j,J}\|f_{j,J}\|_2^2 M_{r_1}(f_{j,J})^{p_b-2}\lesssim
C_* R^\eta
r^{(1-\alpha)\eta+b_0}
\|f\|_2^2M_r(f)^{p_b-2}.
\end{aligned}
\]
By \eqref{eq:step3}, the left-hand side of \eqref{eq:IH} is at most $C(\log r)^{C}r^{-\alpha\eta+b_0}$ times the right-hand side. Since $b_0<\alpha\eta/2$ and $r\ge R_0^\sigma$ is large, this factor is less than one, which closes
\eqref{eq:IH} in this case.

\step{Case 2: The two-ends case.}
Suppose $\beta>r^{b_0}$. We rescale by $x\mapsto x/r$ as at the end of Section~\ref{subsec:two-analytic}, and we use the notation $\delta_*$ and $m_*$ from there. By \eqref{eq:mstar} the tubes of $\T$ form an $m_*$-parallel family of $\delta_*$-tubes, and by \eqref{eq:lambda0} the shading $Y$ is $\lambda_0$-dense. We claim that $Y$ is $(\eps_1,\eps_2)$-two-ends with $\epsilon_1=3\alpha$ and $\epsilon_2=b_0$. Indeed, for large $r$ a ball of radius $\delta_*^{3\alpha}$ corresponds, before the rescaling, to a ball $B(x,\rho)$ with $\rho\le r_1$. Such a ball meets at most three slabs, and $|Y_J(T)|<2\lambda|T|$ for each slab, so by \eqref{eq:lambda0}
\[
 |Y(T)\cap B(x,\rho)|\le6\lambda|T|\le6\beta^{-1}|Y(T)|
 <6r^{-b_0}|Y(T)|\le\delta_*^{b_0}|Y(T)|.
\]
Applying Lemma~\ref{thm:GWX-input} to $Y$, we obtain a set
$E_H\subset X$, which is a union of $r^{1/2+\eps_0}$-cubes, such that $|X\setminus E_H|\lesssim\delta_*^{3\alpha}|X|\le r^{-\alpha}|X|$ and, for every cube $Q\subset E_H$,
\begin{equation}\label{eq:ancestor-H}
 \mu_Y(Q):=\#\{T\in\T:Q\subset Y(T)\}
 \lesssim_{\eps}
 m\lambda_0^{-3/4}
 r^{(n-2)/4+\eps+O(\alpha+\eps_0)} .
\end{equation}
Here we used $m_*\lesssim r^{(n-1)\eps_0}m$ and $\delta_*^{-1}\le r^{1/2}$.

\emph{Removing the cubes outside $E_H$.} Since the quantities $\cB_{\lambda,\beta}(Q)$ are comparable for $Q\in\cQ_2$,
\eqref{eq:refinement-size} gives
\[
\begin{aligned}
 \mathcal B_{\lambda,\beta}(X_2\setminus E_H)
 &\lesssim
 \frac{|X_2\setminus E_H|}{|X_2|}
 \mathcal B_{\lambda,\beta}(X_2)\lesssim r^{-\alpha}(\log r)^{3}
 \mathcal B_{\lambda,\beta}(X_2).
\end{aligned}
\]
For sufficiently large $r$, this is less than half of $\cB_{\lambda,\beta}(X_2)$. We replace $X_2$ by $X_2\cap E_H$ and retain the notation; then \eqref{eq:step3} still holds, and \eqref{eq:ancestor-H} holds for every cube $Q\subset X_2$.

\emph{Fixing a slab, a transverse tuple and the multiplicities.} There are $\lesssim r^{\alpha}$ slabs, so there is an interval $J$ such that the cubes of $X_2$ in the slab of $J$ carry at least the fraction $cr^{-\alpha}$ of $\cB_{\lambda,\beta}(X_2)$. Let $X_J$ be the union of these cubes. By definition, $\cB_{\lambda,\beta}(X_J)$ is the $p$th power of the broad norm on $X_J$, with parameter $A'$, of the single family $f_{\T_J}$. We apply Lemma~\ref{eq:extract-bound} to this family on $X_J$; this is possible because $A'\ge1$. Since there are at most $K^{O(1)}$ transverse tuples, we can fix a transverse tuple of caps
$(\tau_1,\ldots,\tau_k)$ that contributes at least the fraction $K^{-O(1)}$ of the sum. Let $\mathbb T_i$ consist of the packets $T\in\T_J$ with $\tau(T)=\tau_i$, and write
\[
 f_i=\sum_{T\in\mathbb T_i}f_T,
\]
so that $f_i$ is the component of $f_{\T_J}$ with index $\tau_i$. For a cube $Q\subset X_J$ let $\mu_i(Q)$ be the number of tubes of $\T_i$ that meet $Q$. If $\mu_i(Q)=0$ for some $i$, then $\int_Q|Ef_i|^{p}$ is negligible by \eqref{eq:tails}, and so is the contribution of $Q$. For the other cubes, $1\le\mu_i(Q)\le r^{C}$. Pigeonholing the $k$-tuples of dyadic values of these multiplicities, of which there are $O((\log r)^{k})$, we obtain numbers $\mu_1,\ldots,\mu_k\ge1$ and a set
\[
 X_3=\bigcup\{Q\subset X_J:\ \mu_i\le\mu_i(Q)<2\mu_i\ \text{ for }1\le i\le k\}
\]
such that
\begin{equation}\label{boundlow}
\begin{aligned}
 \mathcal B_{\lambda,\beta}(X_2)
 &\lesssim K^{O(1)}r^{O(\alpha)}(\log r)^{O(1)}
 \sum_{b\subset X_3}\min_{1\le i\le k}
 \int_b|Ef_i|^{p_b},
\end{aligned}
\end{equation}
where $b$ runs over the $K^2$-cubes in $X_3\cap B_r$. By H\"older's inequality, as in Lemma~\ref{eq:extract-bound},
\begin{equation}\label{boundhigh}
 \mathcal B_{\lambda,\beta}(X_2)
 \lesssim K^{O(1)}r^{O(\alpha)}(\log r)^{O(1)}
 \left(\prod_{i=1}^k\int_{X_3}|Ef_i|^{p_b}\right)^{1/k}.
\end{equation}

\emph{The quantities attached to $\T_1,\ldots,\T_k$.}
For these families we use the notation of Section~\ref{subsec:two-analytic} and of Lemma~\ref{lem:geometry}, with $X_3$ in the role of $X$ and of $E$: for $T\in\T_i$, $Y_i(T)$ is the union of the cubes of $X_3$ that meet $T$; $N_i=\#\T_i$ and $N=\sum_iN_i$; $\lambda_i$ is the average density of $Y_i$; and $\Lambda$, $\zeta$, $\mu$ are the geometric means defined in Lemma~\ref{lem:geometry}. We record three facts.

First, let $T\in\T_i$. The cubes of $X_3$ belong to $\cQ$ and lie in the slab of $J$, and $T$ is active on $J$; hence every cube of $X_3$ that meets $T$ is contained in $Y_J(T)\subset Y(T)$. Thus $Y_i(T)\subset Y(T)$, and by \eqref{eq:lambda0}
\[
 \lambda_i<4\lambda_0,
 \qquad \Lambda<4\lambda_0.
\]

Second, for the same reason $\mu_i(Q)\le\mu_Y(Q)$ for every cube $Q\subset X_3$. Since $X_3\subset E_H$, the hairbrush estimate \eqref{eq:ancestor-H} gives a bound for each $\mu_i$, in which $\lambda_0^{-3/4}$ may be replaced by $C\Lambda^{-3/4}$:
\begin{equation}\label{eq:mu-hair}
 \frac{\mu_i}{m}\lesssim_\eps\Lambda^{-3/4}r^{(n-2)/4+\eps+O(\alpha+\eps_0)}\qquad(1\le i\le k).
\end{equation}

Third, the sets $\T_i$ are pairwise disjoint subsets of the uniform set $\T^0$. By (W3) and \eqref{eq:am},
\begin{equation}\label{eq:zeta-gain}
 \left(\prod_{i=1}^k\|f_i\|_2^2\right)^{1/k}
 \sim a^2(N_1\cdots N_k)^{1/k}=\zeta a^2N
 \lesssim\zeta\|f\|_2^2.
\end{equation}

\emph{The case of a small multiplicity.}
First suppose $\mu_i<r^{n\eps_0}m$ for some $i$.
Set $t=(p_b-2)/(q-2)$, so that $p_b=2(1-t)+qt$. By \eqref{eq:L2}, \eqref{eq:RD}, $\lambda_i\lesssim1$ and
H\"older's inequality,
\[
 \int_{X_3}|Ef_i|^{p_b}
 \le
 \left(\int_{X_3}|Ef_i|^2\right)^{1-t}
 \left(\int_{X_3}|Ef_i|^q\right)^t\lesssim_{\eps}
 r^{1-2t+\eps+O(\eps_0)}
 \|f_i\|_2^2M_r(f)^{p_b-2}.
\]
Since
\[
 1-2t=-\frac{(d-1)+(n-1)(n-3)}D<0,
\]
and
\[
 \sum_{b\subset X_3}\min_j\int_b|Ef_j|^{p_b}
 \le\int_{X_3}|Ef_i|^{p_b},
\]
\eqref{boundlow}, \eqref{eq:step3} and $\|f_i\|_2^2\lesssim\|f\|_2^2$ give
\eqref{eq:IH}, with the small parameters chosen as above.

\emph{The main case.}
We may therefore assume that
$\mu_i\ge r^{n\eps_0}m$ for every $i$. We verify the hypotheses of Lemma~\ref{lem:geometry} for the rescaled families $\T_1,\ldots,\T_k$, the shadings $Y_i$ and the set $E=X_3$, with $\delta=\delta_*$ and with $m_*$ in place of $m$. Hypothesis~(a) is \eqref{eq:mstar}. For (b), let $T\in\T_i$; its fine cap $\theta(T)$ meets $\supp f_{\tau_i}\subset2\tau_i$, so its center lies in $3\tau_i$, and the direction of $T$ is $G(\xi)$ for some $\xi\in3\tau_i$. By \eqref{eq:transverse}, the wedge condition holds with $\nu=c'_{n,k}K^{-(k-1)}$, and the directions of tubes in different families differ by $\gtrsim K^{-1}\ge\delta_*$. Hypothesis~(c) holds because $\mu_{Y_i}(Q)=\mu_i(Q)\in[\mu_i,2\mu_i)$ for the cubes $Q\subset X_3$, and the densities $\lambda_i$ of Lemma~\ref{lem:geometry} are the ones defined above. Finally $\mu_i/m_*\gtrsim r^{\eps_0}\ge C_{n,k}\log(2/\delta_*)$ for large $r$. Since $Z\lesssim1$, the conclusion \eqref{eq:geom-conclusion} reads
\begin{equation}\label{eq:mu-zahl}
 \frac{\mu}{m}\lesssim_\eps K^{O(1)}\,r^{(n-d)/2+\eps+O(\eps_0)}\,\Lambda^{-(d-1)}\zeta^{-(d-1)} .
\end{equation}

Define
\[
 \mathcal I_\ell
 :=\left(\prod_{i=1}^k\int_{X_3}|Ef_i|^\ell\right)^{1/k},
 \qquad \ell\in\{2,q,p_b\}.
\]
Taking geometric means in \eqref{eq:L2}, in \eqref{eq:RD} combined with \eqref{eq:mu-hair}, and in \eqref{eq:RD} combined with \eqref{eq:mu-zahl}, and using \eqref{eq:zeta-gain} each time, we obtain
\begin{equation}\label{multilinear-norm}
\begin{aligned}
 \mathcal I_2
 &\lesssim
 r^{1+O(\eps_0)}\Lambda\zeta\|f\|_2^2,\\
 \mathcal I_q
 &\lesssim_{\eps}
 r^{-\frac{n}{2(n-1)}+\eps+O(\alpha+\eps_0)}
 \Lambda^{-\frac3{2(n-1)}}\zeta
 \|f\|_2^2M_r(f)^{q-2},\\
 \mathcal I_q
 &\lesssim_{\eps}
 K^{O(1)}r^{-\frac{d-1}{n-1}+\eps+O(\eps_0)}
 \Lambda^{-\frac{2(d-1)}{n-1}}
 \zeta^{1-\frac{2(d-1)}{n-1}}
 \|f\|_2^2M_r(f)^{q-2}.
\end{aligned}
\end{equation}
In the second line, the power of $r$ is $-1+\frac{2}{n-1}\cdot\frac{n-2}{4}=-\frac{n}{2(n-1)}$; in the third, it is $-1+\frac{2}{n-1}\cdot\frac{n-d}{2}=-\frac{d-1}{n-1}$. These are the three bounds discussed in Section~\ref{subsec:strategy}: the $L^2$ estimate, and refined decoupling with the hairbrush estimate and with Zahl's estimate. The factor $\zeta$ in the first two lines is the gain \eqref{eq:zeta-gain}; it compensates for the loss $\zeta^{-2(d-1)/(n-1)}$ in the third line.

Choose
\[
 u=\frac{(2n-3)(d-1)}D,\qquad
 v=\frac{2(n-1)(d-1)}D,\qquad
 w=\frac{(n-1)(n-3)}D.
\]
These weights are positive and sum to one, and they satisfy $2u+q(v+w)=p_b$. A direct computation gives
\[
 u-\frac{n}{2(n-1)}\,v-\frac{d-1}{n-1}\,w=0,
 \qquad
 u-\frac{3}{2(n-1)}\,v-\frac{2(d-1)}{n-1}\,w=0,
\]
so the powers of $r$ and of $\Lambda$ in \eqref{multilinear-norm} cancel, apart from the small losses. Moreover, the resulting power of $\zeta$ is
\[
 1-\frac{2(d-1)}{n-1}w
 =\frac{(2n+1)(d-1)+(n-1)(n-3)}D>0.
\]
By H\"older's inequality, $\mathcal I_{p_b}\le\mathcal I_2^u\mathcal I_q^{v+w}$. Applying the two bounds for $\mathcal I_q$ with weights $v$ and $w$,
respectively, we obtain
\[
\begin{aligned}
 \mathcal I_{p_b}
 &\le\mathcal I_2^u\,\mathcal I_q^{v}\,\mathcal I_q^{w}\lesssim_{\eps}
 K^{O(1)}r^{\eps+O(\alpha+\eps_0)}
 \|f\|_2^2M_r(f)^{p_b-2},
\end{aligned}
\]
since $\zeta\le1$.
Combining this with \eqref{boundhigh} and \eqref{eq:step3} gives
\[
 \|Ef\|_{\BL^{p_b}_{k,A,K}(B_r)}^{p_b}
 \lesssim_{\eps}
 K^{O(1)}r^{\eps+O(\alpha+\eps_0)}(\log r)^{O(1)}
 \|f\|_2^2M_r(f)^{p_b-2}.
\]
By the parameter choices, the right-hand side is at most $R^{\eta/4}r^{\eta/2}\|f\|_2^2M_r(f)^{p_b-2}$ when $r$ is large, which is smaller than the right-hand side of \eqref{eq:IH} if $C_*\ge1$. This proves \eqref{eq:IH} and hence
\eqref{eq:broad-main}.
\end{proof}

\subsection{From mixed broad estimates to linear restriction}
\label{sec:narrow}
Fix $k$ as in \eqref{eq:k-range} and an exponent $p$ with
\begin{equation}\label{eq:p-range}
 \max\{p_b,p_N\}<p<\frac{2(k-1)}{k-2},
\end{equation}
where $p_b=p_b(n,k)$, $p_N=p_N(n,k)$, and the upper endpoint is interpreted as infinity when $k=2$. Such exponents exist. Indeed, $\frac{2(k-1)}{k-2}=2+\frac2{k-2}\ge2+\frac4{n-2}$ by \eqref{eq:k-range}, while $p_N<2+\frac4{n-2}$ because $k<n+2$, and $p_b<2+\frac4{n-2}$ because
\[
 D-(n-2)(2d+n-5)=(2n-1)(d-1)+(n-3)>0 .
\]
We shall prove that for every $\eps>0$ there is a constant $C_\eps$ such that, for all $R\ge1$ and all $f\in L^\infty(B^{n-1})$,
\begin{equation}\label{eq:mixed}
 \int_{B_R}|Ef|^p\le C_\eps R^\eps\|f\|_2^2\|f\|_\infty^{p-2}.
\end{equation}
The proof is by induction on $R$ and follows \cite[Section~9]{Guth}; the only new point is the dependence on the broad parameter $A$, which is now a small power of $R$.

\step{Step 1: the broad estimate at the exponent $p$.} Let $\eta>0$, let $0<\kappa\le\kappa_0(n,k,\eta)$ and $a_0>0$, to be chosen below, and let $K$ be a dyadic number with $R^\kappa\le K<2R^\kappa$ and $A=\lceil R^{a_0}\rceil$. For every cap, $|Ef_\tau|\le\|f_\tau\|_1\lesssim\|f\|_\infty$, and therefore $\int_b|Ef_\tau|^p\lesssim\|f\|_\infty^{p-p_b}\int_b|Ef_\tau|^{p_b}$. Together with Proposition~\ref{prop:broad} and
$M_R(f)\lesssim\|f\|_\infty$, this gives
\begin{equation}\label{eq:broad-p}
 \|Ef\|_{\BL^p_{k,A,K}(B_R)}^p
 \lesssim_{\eta}R^{2\eta}\|f\|_2^2\|f\|_\infty^{p-2}.
\end{equation}

\step{Step 2: the broad--narrow decomposition.} Fix a $K^2$-cube $b\subset B_R$ and subspaces $V_1,\ldots,V_A$ that attain the minimum in the definition \eqref{eq:mu-b} of $\cB_{k,A,p}(f;b)$. We call a cap $\tau$ \emph{broad} if $\ang(G(2\tau),V_a)>K^{-1}$ for all $a$, and \emph{narrow} otherwise. Each narrow cap is assigned to one index $a$ with $\ang(G(2\tau),V_a)\le K^{-1}$; let $\mathcal N_a$ be the set of caps assigned to $a$. Then $Ef=\sum_{\tau\text{ broad}}Ef_\tau+\sum_a\sum_{\tau\in\mathcal N_a}Ef_\tau$, and by H\"older's inequality for the sum over $a$,
\[
 \int_b|Ef|^p\le2^{p-1}\int_b\Bigl|\sum_{\tau\text{ broad}}Ef_\tau\Bigr|^p
 +2^{p-1}A^{p-1}\sum_{a=1}^A\int_b\Bigl|\sum_{\tau\in\mathcal N_a}Ef_\tau\Bigr|^p .
\]
There are $\lesssim K^{n-1}$ caps, so the first term is at most $CK^{(n-1)p}\cB_{k,A,p}(f;b)$. For the second term, the caps of $\mathcal N_a$ lie in the $O(K^{-1})$-neighborhood of a $(k-2)$-dimensional affine subspace of $\R^{n-1}$, and their number is $\lesssim K^{k-2}$. The decoupling inequality of \cite[Lemma~9.3]{Guth}, which is valid for $2\le p\le\frac{2(k-1)}{k-2}$, and H\"older's inequality give, for every $\gamma>0$,
\[
 \int_b\Bigl|\sum_{\tau\in\mathcal N_a}Ef_\tau\Bigr|^p
 \le C_\gamma K^{\gamma+(k-2)(p/2-1)}\sum_{\tau\in\mathcal N_a}\int w_b|Ef_\tau|^p,
\]
where $w_b$ is a weight that is bounded by $1$ and decays rapidly away from $b$. Summing over $a$ and over the cubes $b\subset B_R$, we obtain
\begin{equation}\label{eq:broad-narrow}
\begin{aligned}
 \int_{B_R}|Ef|^p
 \le{}&CK^{(n-1)p}\|Ef\|_{\BL^p_{k,A,K}(B_R)}^p\\
 &+C_{\gamma}A^{p-1}K^{\gamma+(k-2)(p/2-1)}
   \sum_\tau\int_{B_{2R}}|Ef_\tau|^p,
\end{aligned}
\end{equation}
up to a term $\RD(R)\|f\|_2^p$. The factor $A^{p-1}$ comes from
assigning each narrow cap to one of the $A$ subspaces and applying
H\"older's inequality over the resulting groups.

\step{Step 3: parabolic rescaling.} Let $\tau$ be a cap with center $\xi_\tau$ and radius $c_{n,k}K^{-1}$. The change of variables $\xi=\xi_\tau+K^{-1}\tilde\xi$ gives $|Ef_\tau(x)|=K^{-(n-1)}|E\tilde f(\tilde x)|$, where $\tilde f(\tilde\xi)=f_\tau(\xi_\tau+K^{-1}\tilde\xi)$ is supported in $B^{n-1}$, and $\tilde x=(K^{-1}(x'+2x_n\xi_\tau),K^{-2}x_n)$. The map $x\mapsto\tilde x$ has Jacobian $K^{-(n+1)}$ and sends $B_{2R}$ into a ball of radius $CR/K$, and $\|\tilde f\|_2^2=K^{n-1}\|f_\tau\|_2^2$, $\|\tilde f\|_\infty=\|f_\tau\|_\infty$. Covering this ball by $O(1)$ balls of radius $R/K$ and applying \eqref{eq:mixed} at radius $R/K$, we get
\[
 \int_{B_{2R}}|Ef_\tau|^p
 \le C\,C_\eps(R/K)^\eps K^{2n-(n-1)p}
 \|f_\tau\|_2^2\|f_\tau\|_\infty^{p-2}.
\]
Since $\sum_\tau\|f_\tau\|_2^2\|f_\tau\|_\infty^{p-2}\lesssim\|f\|_2^2\|f\|_\infty^{p-2}$ and
\[
 (k-2)\Bigl(\frac p2-1\Bigr)+2n-(n-1)p=-\nu_p,
 \qquad
 \nu_p:=\frac{2n-k}{2}(p-2)-2,
\]
the second term on the right-hand side of \eqref{eq:broad-narrow} is at most
\begin{equation}\label{eq:narrow-term}
 C\,C_\gamma\,C_\eps\,A^{p-1}K^{-\nu_p+\gamma-\eps}R^\eps\|f\|_2^2\|f\|_\infty^{p-2}.
\end{equation}
We have $\nu_p>0$, because $p>p_N=2+\frac4{2n-k}$.

\step{Step 4: closing the induction.} Given $\eps>0$, choose $\eta$, $\kappa\le\kappa_0(n,k,\eta)$, $\gamma$ and $a_0$ such that
\[
 2\eta+\kappa(n-1)p<\frac\eps2,
 \qquad \gamma<\frac{\nu_p}{4},
 \qquad a_0(p-1)<\frac{\kappa\nu_p}{4}.
\]
Then $A^{p-1}\lesssim R^{a_0(p-1)}\le K^{\nu_p/4}$, so the narrow term \eqref{eq:narrow-term} is at most
\[
 C'C_\eps K^{-\nu_p/2}R^\eps\|f\|_2^2\|f\|_\infty^{p-2},
\]
where $C'$ does not depend on $C_\eps$ or $R$. By \eqref{eq:broad-p}, the broad term in \eqref{eq:broad-narrow} is at most $C''R^{\eps/2}\|f\|_2^2\|f\|_\infty^{p-2}$. Since $K\ge R^\kappa$, we have $C'K^{-\nu_p/2}\le\frac12$ for $R\ge R_0$. If $C_\eps\ge2C''$ and $C_\eps$ is large enough for \eqref{eq:mixed} to hold when $R\le R_0$, the induction closes. Note that $a_0$ is chosen after $\kappa$; this is possible because Proposition~\ref{prop:broad} holds for every $a_0>0$.

\step{Step 5: conclusion.} Let $p_1$ be larger than the threshold in Theorem~\ref{thm:main}, let $k$ attain the minimum there, and choose $p<p_1$ in the range \eqref{eq:p-range}. For a measurable set $F\subset B^{n-1}$ and a function $f$ with $|f|\le\1_F$, the estimate \eqref{eq:mixed} gives the local restricted-type bound $\|Ef\|_{L^p(B_R)}\le C_\eps R^{\eps/p}|F|^{1/p}$, and this holds for every exponent in the open interval \eqref{eq:p-range}. The $\eps$-removal argument of Tao \cite{Tao99} converts these local bounds into global restricted-type bounds on $\R^n$ for the exponents in the same interval, since each of them is larger than an exponent for which the local bound holds. Real interpolation between two such exponents then gives the strong-type estimate \eqref{eq:diagonal} for the exponents in the interval, in particular for some exponent below $p_1$, and interpolation with the trivial bound $\|Ef\|_\infty\le\|f\|_1\lesssim\|f\|_\infty$ gives \eqref{eq:diagonal} for $p_1$. This completes the proof of Theorem~\ref{thm:main}.


\begin{thebibliography}{99}

\bibitem{BCT}
J. Bennett, A. Carbery, and T. Tao,
\emph{On the multilinear restriction and Kakeya conjectures},
Acta Math. \textbf{196} (2006), 261--302.

\bibitem{BG}
J. Bourgain and L. Guth,
\emph{Bounds on oscillatory integral operators based on multilinear estimates},
Geom. Funct. Anal. \textbf{21} (2011), 1239--1295.

\bibitem{DW}
C. Demeter and S. Wu,
\emph{Restriction and decoupling estimates for the hyperbolic paraboloid in $\R^3$},
preprint, arXiv:2505.09037v2 (2025).

\bibitem{GWX}
C. Gao, S. Wu, and Y. Xi,
\emph{Sharp microlocal Kakeya--Nikodym estimates for eigenfunctions with applications},
preprint, arXiv:2509.01116v4 (2026).

\bibitem{GWZ}
S. Guo, H. Wang, and R. Zhang,
\emph{A dichotomy for H\"ormander-type oscillatory integral operators},
Invent. Math. \textbf{238} (2024), 503--584.

\bibitem{GuthEndpoint}
L. Guth,
\emph{The endpoint case of the Bennett--Carbery--Tao multilinear Kakeya conjecture},
Acta Math. \textbf{205} (2010), 263--286.

\bibitem{Guth}
L. Guth,
\emph{Restriction estimates using polynomial partitioning II},
Acta Math. \textbf{221} (2018), 81--142.

\bibitem{GIOW}
L. Guth, A. Iosevich, Y. Ou, and H. Wang,
\emph{On Falconer's distance set problem in the plane},
preprint version, arXiv:1808.09346 (2018).

\bibitem{HZ}
J. Hickman and J. Zahl,
\emph{A note on Fourier restriction and nested polynomial Wolff axioms},
J. Anal. Math. \textbf{152} (2024), 19--52.

\bibitem{KT}
N. H. Katz and T. Tao,
\emph{New bounds for Kakeya problems},
J. Anal. Math. \textbf{87} (2002), 231--263.

\bibitem{OAI}
OpenAI,
\emph{Local smoothing for the Euclidean wave equation},
preprint, 2026.
\href{https://github.com/openai/math/blob/main/preprints/Local-smoothing-for-the-Euclidean-wave-equation-October-25-2026/paper.pdf}{Preprint PDF}.

\bibitem{Tao99}
T. Tao,
\emph{The Bochner--Riesz conjecture implies the restriction conjecture},
Duke Math. J. \textbf{96} (1999), 363--376.

\bibitem{WW}
H. Wang and S. Wu,
\emph{Restriction estimates using decoupling theorems and two-ends Furstenberg inequalities},
preprint, arXiv:2411.08871v3 (2024).

\bibitem{Wolff}
T. Wolff,
\emph{An improved bound for Kakeya type maximal functions},
Rev. Mat. Iberoamericana \textbf{11} (1995), 651--674.

\bibitem{Zahl}
J. Zahl,
\emph{New Kakeya estimates using Gromov's algebraic lemma},
Adv. Math. \textbf{380} (2021), 107596.

\end{thebibliography}
\end{document}